\documentclass[letterpaper,10pt, journal]{IEEEtran}

\pdfoutput=1
\usepackage{times} 
\usepackage{amsmath} 
\usepackage{amssymb}  
\usepackage{graphicx,color}
\usepackage{cite}      
\usepackage{array}
\usepackage{times}
\usepackage{latexsym}
\usepackage{amsfonts,amssymb,amsbsy}
\usepackage{mathrsfs}
\usepackage{epsfig} 
\usepackage{enumerate}
\usepackage{framed}
\usepackage{comment}
\usepackage{subfigure}
\definecolor{shadecolor}{rgb}{1.0, 0.92, 0.8}
\usepackage{listings}
\definecolor{Ablue}{rgb}{0.16, 0.32, 0.75}

\newcommand{\mc}{\mathcal}

\newcommand{\mb}{\mathbb}

\newcommand{\diniD}{\mathscr{D}}        

\newcommand{\norm}[1]{\|#1\|}
\def\d{\mathrm{d}}
\def\id{\mathrm{id}}

\newcommand{\C}[1]{\mathbb{#1}}
\newcommand{\R}{\C{R}}

\newtheorem{definition}{Definition}[section]
\newtheorem{theorem}{Theorem}[section]
\newtheorem{remark}{Remark}[section]
\newtheorem{assumption}{Assumption}[section]
\newtheorem{corollary}{Corollary}[section]
\newtheorem{lemma}{Lemma}[section]

\title{\LARGE
Computation of Lyapunov functions for nonlinear differential equations\\
via a Massera--type construction
}
\author{Alina I. Doban, \quad Mircea Lazar  
\thanks{A.I. Doban and  M. Lazar are with the Department of Electrical Engineering, Eindhoven University of Technology, The Netherlands, E-mails: {\tt\small a.i.doban@tue.nl, m.lazar@tue.nl}.
This research is supported
by the Dutch Institute of Systems and Control through a grant from
the NWO graduate program.}%
}

\begin{document}
\maketitle
\thispagestyle{empty}
\pagestyle{empty}
\begin{abstract}
An approach for computing Lyapunov functions for nonlinear continuous--time differential equations is developed via a new, Massera--type construction. This construction is enabled by imposing a finite--time criterion on the integrated function. By means of this approach, we relax the assumptions of exponential stability on the system dynamics, while still allowing integration over a finite time interval. The resulting Lyapunov function can be computed based on any $\mc{K}_\infty$--function of the norm of the solution of the system. In addition, we show how the developed converse theorem can be used to construct an estimate of the domain of attraction. Finally, a range of examples from literature and biological applications such as the genetic toggle switch, the repressilator and the HPA axis are worked out to demonstrate the efficiency and improvement in computations of the proposed approach.
\end{abstract}

\section{Introduction}
\label{sec:01}
The converse of Lyapunov's second method (or direct method) for general nonlinear systems is a topic of extensive ongoing research in the Lyapunov theory community.
Work on the converse theorem started around the $1950$s with the crucial result
in \cite{Massera1949}, which states that if the origin of an autonomous differential equation is asymptotically stable, then the function defined by a semidefinite integral of an appropriately chosen function of the norm of the solution is a continuously differentiable Lyapunov function (LF). This construction led to a significant amount of subsequent work, out of which we recall here \cite{Kurzweil}.
It is well known that finding an explicit form of a LF for general nonlinear systems is a very difficult problem. One of the constructive results on answering the converse problem has been introduced in \cite{Zubov64}, also for autonomous systems. Therein an analytic formula of a LF is provided, which approaches the value $1$ on the boundary of the domain of attraction (DOA) of the considered system. Thus, simultaneously with constructing a LF, an estimate of the DOA is computed. This result, also known as the Zubov method is summarized in \cite[Theorem 34.1]{Hahn67} and \cite[Theorem 51.1]{Hahn67}.
Stemming from Zubov's method, a recursive procedure for constructing a rational LF for nonlinear systems has been proposed in \cite{VanVid85}. This procedure has many computational advantages and it is directly applicable to polynomial systems, providing nonconservative DOA estimates. An alternative construction to the one of Massera, has been proposed  in \cite{Yoshizawa1966},  where it was shown that the supremum of a function of the solutions of the system is a LF.
Additionally, we refer the interested reader to the books \cite{Krasovskii1963} and \cite{LaSalleLef61}, and the survey \cite{Kalman1959}.

As for more recent works, for the particular case of differential inclusions, a converse theorem for uniform global asymptotic stability of a compact set was provided in \cite{LinSontagWang1996}, and for the case of homogeneous systems it was shown in \cite{Rosier92} that asymptotic stability implies the existence of a smooth homogeneous LF. Further converse results for differential inclusions for stability with two measures were provided in, for example, \cite{TeelPraly} and \cite{KelletTeel2004}. If control inputs are to be considered, an existence result of control LFs under the assumption of asymptotic controlability was derived in \cite{Sontag1983}.

For what concerns state--of--the--art, LF constructive methods, see the recent developments of the author of \cite{Hafstein2007} and subsequent works, out of which we single out \cite{Bjornsson2014} and \cite{Bjornsson2015}, where the Massera construction is exploited for generating piecewise affine LFs and \cite{SiggiKelletLi2014} and \cite{Hafstein2016} where the Yoshizawa construction is used.
More detailed historical survey on converse LF results and computatinal methods for LFs can be found in the extensive papers \cite{Kellet2015survey} and \cite{Giesl2015}.

Despite the comprehensive work on the topic of providing a converse to Lyapunov's theorem, the existing constructive approaches either rely on complex candidate LFs (rational, polynomial) or they involve state space partitions (for which scalability with the state space dimension is problematic), accompanied by correspondingly complex or large optimization problems.
In turn, if we restrict strictly to analytical, Massera type of converse results, the construction in \cite[Theorem 4.14]{Khalil2002}, for example, involves integrating over a finite time interval, however with the assumption of exponential stability of the origin. A similar construction, with the relaxed assumption of asymptotic stability, has been developed in \cite{Bjornsson2015}, by using a Lipschitz, positive outside a neighborhood around the origin, (arbitrary) function of the state.

In this paper, we propose a similar Massera--type construction for the LF, by relaxing the exponential stability  assumption to a richer type of $\mc{K}\mc{L}$--stability property. Additionally, we allow for the LF to be generated by any  $\mc{K}_\infty$ candidate function which satisfies a
finite--time decrease criterion. Nonetheless, this relaxation comes with a restriction on the $\mc{K}\mc{L}$--stability condition indicated in the paper by Assumption~\ref{as:02.01}. Ultimately, the proposed solution makes use of an analytic relation between LFs and finite--time Lyapunov functions (FTLFs) to compute a LF. Thus, construction of LFs is brought down to verification that a candidate function is a FTLF, which is somewhat easier than identifying a candidate function for a true LF.

A similar finite--time criterion has been introduced in \cite{Aeyels1998} to provide a new asymptotic stability result for nonautonomous nonlinear differential equations. The discrete--time analog of this condition has been first used in \cite{Geiselhart2014} to provide a converse Lyapunov theorem for nonlinear difference equations. The candidate LF therein is also of Massera type, but projected in discrete--time, thus defined by finite summation. In this respect, some of the results reported in this paper provide a continuous--time counterpart to the findings in \cite{Geiselhart2014}.

The main results of this paper consist of the finite--time converse result in Theorem~\ref{th:02.02}, the equivalence condition in Lemma~\ref{lemma:02.01} and the construction in Theorems~\ref{co:03.01} and \ref{th:04.01}. One of the major benefits of the proposed converse theorem is that it enables a systematic construction of DOA estimates.
Since the computation procedure is based on analytical relations, it provides potentially improved scalability with the state--space dimension and it is applicable also to systems with nonpolynomial nonlinearities.

\subsection{Notation and some definitions}
\label{sec:01.01}
We say that a set $\mc{S}\subseteq\R^n$ is  proper if it contains the origin in its interior and it is compact. The logarithmic norm \cite{Soderlind} of a matrix $A\in\R^{n\times n}$,  $\mu(A)$\footnote{Thus, $\mu(A)$, sometimes reffered to as the matrix measure, does not define a norm in the conventional sense.} is defined as $\mu(A)=\underset{h\rightarrow 0^+}\lim\frac{\norm{I+hA}-1}{h}$. For the $1,2$ and $\infty$ norms, standard definitions of $\mu(A)$ exist. We recall here the definition of the logarithmic norm induced by the $2$--norm, $\mu_2(A)$ and by the $2$--weighted norm, $\mu_{2,P}(A)$ \cite{Hu2004}:
$\mu_2(A)=\lambda_{max}(\frac{1}{2}(A+A^\top))$, $\mu_{2,P}(A)=\lambda_{max}(\frac{\sqrt{P}A\sqrt{P}^{-1} + (\sqrt{P}A\sqrt{P}^{-1})^\top}{2})$, where $\lambda_{max}(P)$ denotes the largest eigenvalue of a symmetric real matrix $P$.

In this paper we consider autonomous continuous--time systems described by
\begin{equation}
\label{eq:01.01}
\dot{x}=f(x),
\end{equation}
where $f\,:\,\R^n\,\rightarrow\,\R^n$, is a locally Lipschitz function.
%
\begin{remark}[Solution notation.]
\label{ch4:re:sol}
Let the solution of \eqref{eq:01.01} at time $t\in\R_{\geq 0}$ with initial value $x(0)$ be denoted by $\phi(t,x(0))$, where $\phi\,:\,\R_{\geq 0}\times\R^n\,\rightarrow \R^n$. 
We assume that $\phi(t,x(0))$ exists and it is unique for all $t\in\R_{\geq0}$ (see \cite[Chapter 3]{Khalil2002} for sufficient smoothness conditions on $f$).
The locally Lipschitz assumption on $f(x)$ implies, additionally,  that $\phi(t,x(0))$ is a continuous function of $x(0)$ \cite[Chapter III]{Hahn67}.
Furthermore, we assume that the system \eqref{eq:01.01} has an equilibrium point at the origin, i.e. $f(0)=0$. 

In what follows, for simplicity, unless stated otherwise, we will use the notation $x(t):=\phi(t,x(0))$ with the following interpretations:
\begin{itemize}
\item if $x(t)$ is the argument of $\dot{W}(x(t))$, then $x(t)$ represents the solution of the system as defined by $\phi(t, x(0))$;
\item if $x(t)$ is the argument of $V(x(t))$ as in \eqref{eq:02.01b} and \eqref{eq:02.06}, for example,  then $x(t)$ represents a point on the solution $\phi(t,x(0))$ for a fixed value of $t$; the same interpretation holds for $W(x(t))$; in other words we do not consider time-varying functions $V$ or $W$.
\end{itemize}
\end{remark}
In what follows, we proceed by recalling some subsidiary notions and definitions.
\begin{definition}
\label{def:01.01}
A function $\alpha\,:\,\R_{\geq0}\,\rightarrow\R_{\geq 0}$ is said to be a $\mc{K}$ function if it is continuous, zero at zero and strictly increasing.
If additionally, $\lim_{s\rightarrow\infty}\alpha(s)=\infty$, then $\alpha$ is called a $\mc{K}_\infty$ function.
\end{definition}
\begin{definition}
\label{def:01.02}
A function $\sigma\,:\,\R_{\geq0}\,\rightarrow\R_{\geq 0}$ is said to be a $\mc{L}$ function  if it is continuous, strictly decreasing and  $\lim_{s\rightarrow\infty}\sigma(s)=0$.
\end{definition}
\begin{definition}
\label{def:01.03}
A function $\beta\,:\,\R_{\geq0}\times \R_{\geq0}\,\rightarrow\R_{\geq 0}$ is said to be a $\mc{K}\mc{L}$ function  if it is a $\mc{K}$ function in its first argument and a $\mc{L}$ function in its second argument.
\end{definition}
\begin{definition}
\label{def:01.04}
The origin is an asymptotically stable (AS) equilibrium for the system \eqref{eq:01.01} if for some proper set $\mc{S}\subseteq\R^n$, there exists a function $\beta\in\mc{K}\mc{L}$  such that for all $x(0)\in\mc{S}$,
\begin{equation}
\label{eq:01.02}
\norm{x(t)}\leq\beta(\norm{x(0)},t),\quad \forall t\in\R_{\geq 0}.
\end{equation}
If the set $\mc{S}=\R^n$, then we say that the origin is globally asymptotically stable (GAS).
\end{definition}

AS defined as above is equivalent with $\mc{K}\mc{L}$--stability in the set $\mc{S}$. In the remainder of the paper we will say that the origin is \emph{$\mc{K}\mc{L}$--stable in $\mc{S}$} to refer to the property defined above. When $\mc{S}=\R^n$, then we use the term global $\mc{K}\mc{L}$--stability.
\begin{definition}
\label{def:01.05}
A continuously differentiable function $V\,:\,\R^n\rightarrow\R_{\geq0}$, for which there exist $\alpha_1,\alpha_2\in\mc{K}_\infty$ and a $\mc{K}$ function $\rho \,:\,\R_{\geq 0}\rightarrow\R_{\geq0}$ such that
\begin{eqnarray}
\label{eq:01.03a}
\alpha_1(\norm{x}) \leq  V(x) &\leq& \alpha_2(\norm{x}), \quad \forall x\in\R^n\\
\label{eq:01.03b}
                \dot{V}(x)=\nabla^\top V(x) f(x)&\leq&-\rho(\norm{x}), \quad \forall x\in\mc{S},
\end{eqnarray}
with $\mc{S}\subseteq\R^n$ proper, is called a Lyapunov function for the system \eqref{eq:01.01}.
\end{definition}
\begin{definition}
\label{def:01.05}
A proper set $\mc{S}\subseteq\R^n$ is called an invariant set for the system \eqref{eq:01.01} if for any $x(0)\in\mc{S}$, the corresponding solution $x(t)\in\mc{S}$, for all $t\in\R_{\geq 0}$.
\end{definition}
\begin{definition}
\label{def:02.01}
Given a positive, real scalar $d$, the proper set $\mc{S}\subseteq\R^n$ is called a $d$--invariant set for the system \eqref{eq:01.01} if for any  $t\in \R_{\geq 0}$, if $x(t)\in\mc{S}$, then it holds that $x(t+d)\in\mc{S}$.
\end{definition}
Note that the $d$--invariance property does not imply that $x(t)\in\mc{S}$ for all $t\geq 0$ if $x(0)\in\mc{S}$.

We recall below Sontag's lemma on $\mc{K}\mc{L}$--estimates \cite[Proposition 7]{Sontag98}, as it will be instrumental.
\begin{lemma}
\label{lemma:01.01}
For each class $\mc{K}\mc{L}$--function $\beta$ and each number $\lambda\in\R_{\geq 0}$, there exist $\varphi_1, \varphi_2\in\mc{K}_\infty$, such that $\varphi_1(s)$ is locally Lipschitz and
\begin{equation}
\label{eq:01.04}
\varphi_1(\beta(s,t))\leq\varphi_2(s)e^{-\lambda t}, \quad\forall s, t\in\R_{\geq 0}.
\end{equation}
\end{lemma}

The following result was introduced in \cite[Definition 24.3]{Hahn67} to relate positive definite functions and  $\mc{K}$--functions. 
A proof was proposed in \cite[Lemma 4.3]{Khalil2002}.

\begin{lemma}
\label{lemma:01.02}
Consider a function $W\,:\,\R^n\,\rightarrow\,\R_{\geq0}$ with $W(0)=0$.
\begin{enumerate}[1.]
\item
If $W(x)$ is continuous and positive definite in some neighborhood around the origin, $\mc{N}(0)$, then there exist two functions
$\hat{\alpha}_1, \hat{\alpha}_2\in\mc{K}$  such that
\begin{equation}
\label{eq:01.05}
\hat{\alpha}_1(\norm{x})\leq W(x)\leq\hat{\alpha}_2(\norm{x}), \quad \forall x\in\mc{N}(0).
\end{equation}
\item
If $W(x)$ is continuous and positive definite in $\R^n$ and additionally, $W(x)\rightarrow\infty$, when $x\rightarrow\infty$ then \eqref{eq:01.05} holds with $\hat{\alpha}_1, \hat{\alpha}_2\in\mc{K}_\infty$ and for all $x\in\R^n$.
\end{enumerate}
\end{lemma}
Next we recall the Bellman--Gronwall Lemma. A proof is provided in \cite[Lemma C. 3.1]{Sontag1998}.
\begin{lemma}
\label{lemma:01.03}
Assume given an interval $\mc{I}\subseteq\R$, a constant $c\in\R_{\geq 0}$, and two functions $\alpha,\mu\,:\,\mc{I}\,\rightarrow\,\R_{\geq 0}$, such that $\alpha$ is locally integrable and $\mu$ is continuous. Suppose further that for  some $\sigma\in\mc{I}$ it holds that
$$\mu(t)\leq\nu(t):=c+\int_\sigma^t \alpha(\tau)\mu(\tau)\d\tau$$
for all $t\geq \sigma$, $t\in\mc{I}$. Then it must hold that
$$\mu(t)\leq ce^{\int_\sigma^t\alpha(\tau)\d\tau}.$$
\end{lemma}

It is well known, from the direct method of Lyapunov, that the existence of a Lyapunov function for the system \eqref{eq:01.01} implies that the origin is an AS equilibrium for \eqref{eq:01.01}. In the remainder of this paper we propose some new alternatives to classical Lyapunov converse results, which are verifiable and constructive towards obtaining nonconservative DOA estimates. More specifically, in Section~\ref{sec:02.01} we recall a finite--time criterion for $\mc{K}\mc{L}$--stability in a given set $\mc{S}$ and we provide its converse. In Section~\ref{sec:02.02} the alternative converse theorem is provided and comparative remarks with respect to existing constructions are drawn in Section~\ref{sec:02.03}.
Given a FTLF $V$, an expansion scheme which preserves the FT decrease property is rendered in Section~\ref{sec:03}, while indirect verification techniques are given in Section~\ref{sec:04.01} and computational steps are indicated in Section~\ref{sec:04.02}. Finally, a range of insightful worked out examples from literature and biological applications (biochemical reactions) such as the genetic toggle switch and the HPA axis are shown in Section~\ref{sec:05}.

\section{A constructive Lyapunov converse theorem}
\label{sec:02}
\subsection{Finite--time conditions}
\label{sec:02.01}
\begin{definition}
\label{def:finite}
Let there be a continuous function $V\,:\,\R^n\,\rightarrow\,\R_{\geq 0}$, and a real scalar $d>0$  for which the proper set $\mc{S}\subseteq\R^n$ is $d$--invariant and
the conditions
\begin{eqnarray}
\label{eq:02.01a}
\alpha_1(\norm{x})\leq V(x)&\leq &\alpha_2(\norm{x}), \quad\forall x\in\R^n,\\
\label{eq:02.01b}
V(x(t+d))-V(x(t))&\leq& -\gamma(\|x(t)\|) ,\quad \forall t\geq 0,
\end{eqnarray}
are satisfied with $\alpha_1$, $\alpha_2\in\mc{K}_\infty$, and $\gamma\in\mc{K}$ and for all $x(t)$, with $x(0)\in\mc{S}$.
Then the function $V$  is called a \emph{finite--time Lyapunov function} (FTLF) for the system \eqref{eq:01.01}.
\end{definition}

In order for condition \eqref{eq:02.01b} to be well--defined, additionally to the locally Lipschitz property of the map $f(x)$, it is assumed that there exists no finite escape time in each interval $[t,t+d]$, for all $t\in\R_{\geq 0}$. However, as it will be shown later, it is sufficient to require that there is no finite escape time in the time interval $[0,d]$.

When $\mc{S}=\R^n$, a sufficient condition for existence of the solution for all $t\in\R_{\geq 0}$ is that the map $f(x)$ is Lipschitz bounded \cite[Chapter III.16]{Hahn67}. Furthermore, note that existence of a finite escape time for initial conditions in a given set in $\R^n$ implies that the origin is unstable in that set \cite[Chapter III.16]{Hahn67}.

The following result relates inequality \eqref{eq:02.01b} with another known type of decrease condition, which will be instrumental.


\begin{lemma}
\label{re:02.01}
The decrease condition \eqref{eq:02.01b} on $V$ is equivalent with
\begin{equation}
\label{eq:02.06}
V(x(t+d))-\rho(V(x(t)))\leq 0,\quad \forall t\in\R_{\geq 0},
\end{equation}
for all $x(t)$ with $x(0)\in\mc{S}$, where $\rho\,:\,\R_{\geq 0}\,\rightarrow\,\R_{\geq 0}$ is a positive definite, continuous function and satisfies $\rho<\id$ and $\rho(0)=0$.
\end{lemma}
\begin{IEEEproof}
The proof follows a similar reasoning as in \cite[Remark 2.5]{Geiselhart2015}.
Assume that $V$ is such that \eqref{eq:02.01a} and \eqref{eq:02.01b} hold. Then, for all $x(t)\neq0$:
\begin{equation}
\label{eq:02.07}
\begin{split}
0\leq V(x(t+d))&\leq V(x(t))-\gamma(\norm{x(t)})\\
               &< V(x(t))-0.5\gamma(\norm{x(t)})\\
               &\leq V(x(t))-0.5\gamma(\alpha_2^{-1}(V(x(t))))\\
               &=(\id-0.5\gamma\circ \alpha_2^{-1})(V(x(t)))\\
               &=:\rho(V(x(t))).
\end{split}
\end{equation}
Similarly, for all $x(t)\neq 0$, 
\begin{equation}
\nonumber
\begin{split}
0\leq V(x(t+d))&\leq V(x(t))-\gamma(\norm{x(t)})\\
               &< \alpha_2(\norm{x(t)})-0.5\gamma(\norm{x(t)})\\
               &=(\alpha_2-0.5\gamma)(\norm{x(t)}),\\
\end{split}
\end{equation}
which implies that $(\alpha_2-0.5\gamma)(s)>0$, for all $s\neq0$. Furthermore, since $\alpha_2^{-1}\in\mc{K}_\infty$, then $(\alpha_2-0.5\gamma)\circ\alpha_2^{-1}(s)>0$ and
$$0<(\id-0.5\gamma\circ\alpha_2^{-1})(s)<\id, \quad \forall s\neq0.$$
Thus, by construction, the function $\rho\,:\,\R_{\geq 0}\,\rightarrow\,\R_{\geq 0}$ is a continuous, positive definite function. All involved functions are continuous by definition, thus the difference remains continuous. When $x(t)=0$, then \eqref{eq:02.06} trivially holds, since $\rho(0)=0$.  

Now assume that \eqref{eq:02.06} holds. Then
\begin{equation}
\label{eq:02.08}
\begin{split}
V(x(t+d))-V(x(t))&\leq \rho(V(x(t)))-V(x(t))\\
                 &=-(V(x(t))-\rho(V(x(t)))\\
                 &=-((\id-\rho)(V(x(t))))\\
                 &\leq-((\id-\rho)(\alpha_1(\norm{x(t)}))\\
                 &=-\tilde{\gamma}(\norm{x(t)}),
\end{split}
\end{equation}
with $\tilde{\gamma} = (\id-\rho)\circ\alpha_1$. Since $\rho<\id$ by assumption, then $\tilde{\gamma}$ is positive definite, and furthermore continuous. Thus, by Lemma~\ref{lemma:01.02} $\tilde{\gamma}$ can be can be lower bounded by a $\mc{K}$--function $\gamma(\norm{x(t)})$,  hence
$V(x(t+d))-V(x(t))\leq -\gamma(\norm{x(t)})$, with $\gamma\in\mc{K}.$
\end{IEEEproof}

Next, we propose a version of \cite[Theorem $1$]{Aeyels1998} for the time--invariant case and with the additional assumption that the set $\mc{S}$ is a $d$--invariant set for \eqref{eq:01.01}.
This additional assumption enables a simpler proof, while the result is stronger, i.e., $\mc{K}\mc{L}$--stability in $\mc{S}$ is attained as opposed to local (in some neighborhood around the origin) $\mc{K}\mc{L}$--stability.
\begin{theorem}
\label{th:02.01a}
If a function $V$ defined as in \eqref{eq:02.01a} and \eqref{eq:02.01b} and a proper $d$-invariant set $\mc{S}$ exist for the system \eqref{eq:01.01}, then the origin equilibrium of \eqref{eq:01.01} is $\mc{K}\mc{L}$--stable in $\mc{S}$.
\end{theorem}
\begin{IEEEproof}
For any $t\in\R_{\geq 0}$, there exists an integer $N\geq0$ and $j\in\R_{\geq0}$, $j<d$ such that $t=Nd+j$.
By applying \eqref{eq:02.01b} in its equivalent form \eqref{eq:02.06} recursively, we get that
\begin{equation}
\label{eq:02.AS}
\begin{split}
V(x(t))&=   V(x(Nd+j))\\
       &=   V(x(((N-1)d+j)+d))\\
       &\leq\rho( V(x((N-1)d+j)))\\
       &=  \rho( V(x(((N-2)d+j)+d)))\\
       &\leq \rho^2(V(x((N-2)d+j)))\\
       &\ldots\\
       &\leq \rho^N(V(x(j)))\\
       &\leq\rho^N(\alpha_2(\norm{x(j)})),
\end{split}
\end{equation}
where $\rho^N$ denotes the $N$--times composition of $\rho$.
The solution at time $t=j$ is given by
$$x(j) = x(0)+\int_{0}^j f(x(s))\d s, $$
for any $j\geq 0$.
Then
\begin{equation*}
\begin{split}
\norm{x(j)-x(0)}\leq &\int_0^j \norm{f(x(s))-f(x(0))+f(x(0))}\d s \\
                \leq &\int_0^j \norm{f(x(s))-f(x(0))}\d s +\\
                &\int_0^j\norm{f(x(0))}\d s.
\end{split}
\end{equation*}
By using the local Lipschitz continuity property of $f$, with $L>0$ the Lipschitz constant,  and the Lemma~\ref{lemma:01.03} we obtain that
\begin{equation*}
\begin{split}
\norm{x(j)-x(0)}& \leq \int_0^j L\norm{x(s)-x(0)}\d s +\int_0^j\norm{f(x(0))}\d s\\
                &\leq  \left(\int_0^j\norm{f(x(0))}\d s\right)e^{Lj}.
\end{split}
\end{equation*}
Thus,
$$\norm{x(j)}\leq\norm{x(0)}+\left(\int_0^j\norm{f(x(0))}\d s\right)e^{Lj}=:F_j(\norm{x(0)}).$$
Then $F_j(\norm{x(0)})\leq F_d(\norm{x(0)})$, for all $j\in[0,d]$.
By the standing assumptions on $f$ it results that $F_d(\norm{x(0)})$ is continuous with respect to $x(0)$.
Furthermore, $F_d(0)=0$, $F_d(s)$ is positive definite and continuous and $F_d(s)\rightarrow\infty$, when $s\rightarrow\infty$, for any $s\geq 0$. By applying Lemma~\ref{lemma:01.02} to $F_d(\norm{x(0)})$ we obtain that there exists a function $\omega\in\mc{K}_\infty$ such that $F_d(\norm{x(0)})\leq\omega(\norm{x(0)})$, and consequently, 
$\norm{x(j)}\leq\omega(\norm{x(0)})$, for all $0\leq j<d$.
Thus, with $\hat{\alpha}_2:=\alpha_2\circ\omega$ and $\hat{\rho}:=\rho^{-1}$ we get that
\begin{equation}
\nonumber
\begin{split}
V(x(t))&\leq\rho^N(\hat{\alpha}_2(\norm{x(0)}))
       =   \rho^{\frac{t-j}{d}}(\hat{\alpha}_2(\norm{x(0)}))\\
       &\leq\rho^{\lfloor \frac{t}{d}\rfloor -1}\circ\hat{\alpha}_2(\norm{x(0)})\\
       &=\rho^{\lfloor \frac{t}{d}\rfloor}\circ\rho^{-1}\circ\hat{\alpha}_2(\norm{x(0)})\\
       &\leq\rho^{\lfloor \frac{t}{d}\rfloor}\circ\hat{\rho}\circ\hat{\alpha}_2(\norm{x(0)}),\quad \hat{\rho}\in\mc{K}_\infty\\
       &=:\hat{\beta}(\norm{x(0)},t).
\end{split}
\end{equation}
Without loss of generality we can assume that $\rho$ is a one--to--one (injective) and onto (surjective) function, thus invertible since, by hypothesys, $\mc{S}$ is  compact. Furthermore, since $\rho$ is continuous, then by \cite[Theorem 3.16]{Browder1996}, $\rho^{-1}$ is continuous. Additionally, $\rho^{-1}(0)=\rho^{-1}(\rho(0))=0$. Thus, there exists a function $\hat{\rho}\in\mc{K}_{\infty}$, such that $\rho^{-1}\leq\hat{\rho}$, as follows from Lemma~\ref{lemma:01.02}. We can conclude that $\hat{\beta}\in\mc{K}\mc{L}$ since $\hat{\rho}\circ\hat{\alpha}_2(s)\in\mc{K}_\infty$ and $\rho^{\lfloor \frac{t}{d}\rfloor}\in\mc{L}$.

Finally,
$$\norm{x(t)}\leq\alpha_1^{-1}(\hat{\beta}(\norm{x(0)},t))=:\beta(\norm{x(0)},t),$$
for all $x(0)\in\mc{S}$ and for all $t\in\R_{\geq 0}$,
thus we have obtained $\mc{K}\mc{L}$--stability in $\mc{S}$.
\end{IEEEproof}

A similar result was derived in \cite[Proposition 2.3]{Karafyllis2012}, which  offers an alternative to the periodic decrease condition in \cite{Aeyels1998} (here \eqref{eq:02.01b}), by requiring the minimum over a finite time interval of a positive definite function of the state to decrease. Condition (2.2) in \cite{Karafyllis2012} always implies a decrease after a finite time interval, but it allows the length of the time interval to be state dependent. Proposition 2.3 of \cite{Karafyllis2012} shows that such a relaxed finite time decrease condition implies $\mc{K}\mc{L}$-stability and exponential stability (under the usual global exponential stability assumptions plus a common time interval length for all states).

We proceed by providing a converse finite--time Lyapunov function for $\mc{K}\mc{L}$--stability in a compact set $\mc{S}$.
\subsection{Alternative converse theorem}
\label{sec:02.02}
\begin{assumption}
\label{as:02.01}
There exists a $\mc{K}\mc{L}$--function $\beta$ satisfying \eqref{eq:01.02} for the system \eqref{eq:01.01} such that
\begin{equation}
\label{eq:02.08}
\beta(s, d)<s
\end{equation}
for some positive $d\in\R$ and all $s>0$.
\end{assumption}
\begin{theorem}
\label{th:02.02}
If the origin is $\mc{K}\mc{L}$--stable in some invariant subset of $\R^n$, $\mc{S}$\footnote{Invariance is needed in order for \eqref{eq:02.06} to hold for all $t\geq 0$.} for the system \eqref{eq:01.01} and Assumption~\ref{as:02.01} is satisfied, then for any function $\eta\in\mc{K}_\infty$ and for any norm $\norm{\cdot}$, the function $V:\R^n\,\rightarrow\,\R_{\geq0}$, with
\begin{equation}
\label{eq:02.10}
V(x):=\eta(\norm{x}), \quad \forall x\in\R^n
\end{equation}
satisfies  \eqref{eq:02.01a} and \eqref{eq:02.01b}.
\end{theorem}
\begin{IEEEproof}
Let the pair $(\beta, d)$ be such that Assumption~\ref{as:02.01} holds. Then, by hypothesis we have that:
\begin{equation}
\nonumber
\begin{split}
\eta(\norm{x(t+d)})  
                   &\leq\eta(\beta(\norm{x(t)},d))\\
                   &\leq\eta(\beta(\eta^{-1}(V(x(t))),d))\\
                   &:=\rho(V(x(t))),
\end{split}
\end{equation}
where $\rho =\eta(\beta(\eta^{-1}(\cdot),d))$, for all initial conditions $x(0)\in\mc{S}$. By Assumption~\ref{as:02.01}, we obtain that there exists a $d>0$ such that $\rho<\eta(\eta^{-1}(\cdot))=\id$, Thus, we get
$$V(x(t+d))-\rho(V(x(t)))\leq0, \quad \forall x(0)\in\mc{S}.$$
From Lemma~\ref{re:02.01} this implies that \eqref{eq:02.01b} holds.
Since $V$ is defined by a $\mc{K}_\infty$ function, then let $\alpha_1(s)= \alpha_2(s)=\eta(s)$ such that \eqref{eq:02.01a} holds.
\end{IEEEproof}
In \cite[Remark 2.4]{Karafyllis2012},  a converse result for $\mc{K}\mc{L}$--stable systems is derived, in terms of the finite decrease condition (2.2) therein. More precisely, it is shown that if the $\mc{K}\mc{L}$--stability property holds, then any positive definite function satisfies inequality (2.2) in \cite{Karafyllis2012}. Compared to the converse results of \cite{Karafyllis2012}, the converse theorem above shows that a stronger condition holds (inequality \eqref{eq:02.01b} with a common finite time $d$ for all states $x(t)$) under Assumption~\ref{as:02.01}.

Consider the function defined as
\begin{equation}
\label{eq:02.11}
 W(x(t)):=\int_{t}^{t+d}V(x(\tau))\d\tau,
\end{equation}
for any $V$ that satisfies \eqref{eq:02.01a} and \eqref{eq:02.01b}.

Generally, in standard converse theorems the function $\varphi_1$ which defines the LF is a particular, special $\mc{K}_\infty$ function.
In the proposed finite--time converse theorem, $\eta$ is any $\mc{K}_\infty$ function, which allows for more freedom in the construction.
In turn, the developed finite--time converse theorem will be used to obtain an alternative converse Lyapunov theorem.
\begin{lemma}
\label{lemma:02.01}
A continuously differentiable function $V\,:\,\R^n\,\rightarrow\,\R_{\geq 0}$ satisfies \eqref{eq:02.01a} and \eqref{eq:02.01b} for \eqref{eq:01.01} and some $d>0$, if and only if the function $W$  as defined in \eqref{eq:02.11} with the same $d>0$  is a Lyapunov function for the system \eqref{eq:01.01}.
\end{lemma}
\begin{IEEEproof}
Let there be a function $V$ satisfying \eqref{eq:02.01a} and \eqref{eq:02.01b}.
$V$ is continuous, thus it is integrable over any closed, bounded interval
$[t, t+d]$, $t\geq0$. By Theorem~$5.30$  in \cite{Browder1996}, this implies that $W(x(t))$ is continuous on each  interval $[t, t+d]$, for any $t$.
Since $V$  is also positive definite, by integrating over the bounded interval $[t, t+d]$ the resulting function $W(x(t))$ will also be positive definite.

Additionally, $\lim_{x\rightarrow\infty}V(x)=\infty$, thus,  $\lim_{x\rightarrow\infty}W(x)=\infty$ and 
the result in Lemma~\ref{lemma:01.02} can be applied.
 Therefore, there exist two functions $\hat{\alpha}_1, \hat{\alpha}_2\in\mc{K}_\infty$ such that
\begin{equation}
\label{eq:02.12}
\hat{\alpha}_1(\norm{x})\leq W(x)\leq\hat{\alpha}_2(\norm{x}), \quad \forall x\in\mb\R^n,
\end{equation}
holds.
Next, by making use of the general Leibniz integral rule, we get that
\begin{equation}
\nonumber
\begin{split}
\frac{\d }{ \d t}W(x(t))=&\int_{t}^{t+d}\underbrace{\frac{\d}{\d t} V(x(\tau))}_{=0}\d \tau +\\
                       & V(x(t+d))\dot{(t+d)}-V(x(t))\dot{t}\\
                       =&V(x(t+d))-V(x(t))\leq-\gamma(\norm{x(t)}).
\end{split}
\end{equation}
Thus, $W$ is a Lyapunov function for \eqref{eq:01.01}.

Now assume that $W$ is a Lyapunov function for \eqref{eq:01.01}, i.e.  \eqref{eq:02.12} holds and for $\gamma\in\mc{K}$ it holds that
 $$\dot{W}(x)\leq -\gamma(\norm{x}),\quad \forall x\in\mc{S}.$$
By the same Leibniz rule, we know that $\dot{W}=V(x(t+d))-V(x(t))$, thus the difference $V(x(t+d))-V(x(t))$ is negative definite, i.e. \eqref{eq:02.01b} holds. Now we have to show that  \eqref{eq:02.01a} holds.

%
%

Assume that there exists an $x(t)\in\mc{S}$, $x(t)\neq 0$,  such that $V(x(t))\leq0$. Then, this  implies that
$$V(x(t+d))<V(x(t)\leq0,$$ and furthermore, 
$V(x(t+id))<V(x(t+(i-1)d))<\ldots<V(x(t+2d))<V(x(t+d))<V(x(t))\leq 0,$
for all integers $i>0$, due to $d$-invariance of $\mc{S}$ and the assumption that $x(t) \in \mc{S}$. Then, 
$$\lim_{i\rightarrow\infty}V(x(t+id))=-\infty.$$
Since $W$ is a LF for \eqref{eq:01.01} in $\mc{S}$, then the origin is  $\mc{K}\mc{L}$--stable in $\mc{S}$, thus it implies that $\lim_{i\rightarrow\infty}x(t+id)=0$. Then, because the solution of the system \eqref{eq:01.01} is a continuous function of time and $W$ is continuous, it follows that
$\lim_{i\rightarrow\infty}W(x(t+id))=W(0)$. Therefore, we have that
\begin{equation*}
\begin{split}
\lim_{i\rightarrow\infty}W(x(t+id))&=\lim_{i\rightarrow\infty}\int_{t+id}^{t+(i+1)d}V(x(\tau))\d\tau\\
&\Leftrightarrow\\
W(0) &= \int_\infty^\infty V(x(\tau))\d\tau=V(x(\infty))=-\infty,\\
\end{split}
\end{equation*}
 which is a contradiction since $W(0)=0$,
thus $V(x)$ must be positive definite on $\mc{S}$. 
 By the definition of $W$, we have that $V$ must be a continuous function, because it needs to be integrable for $W$ to exist. By assumption, $W$ is upper and lower bounded by $\mc{K}_\infty$ functions, thus for $x\rightarrow\infty$, $W(x)\rightarrow\infty$. This can only happen when $V(x)\rightarrow\infty$.
Thus, using a similar reasoning as above, based on  Lemma~\ref{lemma:01.02}, this implies that $V$ is upper and lower bounded by $\mc{K}_\infty$ functions, hence \eqref{eq:02.01a} holds.
\end{IEEEproof}

The next result summarizes the proposed alternative converse LF for $\mc{K}\mc{L}$--stability in $\mc{S}$, enabled by the finite--time conditions \eqref{eq:02.01a} and \eqref{eq:02.01b}.
\begin{corollary}
\label{co:02.01}
If the origin is $\mc{K}\mc{L}$--stable in some set $\mc{S}$ for the system \eqref{eq:01.01}, with the $\mc{K}\mc{L}$--function $\beta$ satisfying Assumption~\ref{as:02.01} for some $d>0$, then by Theorem~\ref{th:02.02} and Lemma~\ref{lemma:02.01}, for any function $\eta\in\mc{K}_\infty$ and any norm $\norm{\cdot}$, the function $W(\cdot)$ defined as in \eqref{eq:02.11} for the same $d>0$, is a Lyapunov function for the system \eqref{eq:01.01}. 
\end{corollary}

The above corollary provides a continuous--time counterpart of the converse result
\cite[Corollary 22]{Geiselhart2014}. The main line of reasoning relies on the Assumption~\ref{as:02.01}, which is the same assumption needed for the discrete--time converse theorem.  The main technical differences with \cite{Geiselhart2014} lie in the proof of Theorem~\ref{th:02.01a}, the construction \eqref{eq:02.11} and the proof of Lemma~\ref{lemma:02.01}.

\subsection{Remarks Massera construction}
\label{sec:02.03}
Notice that the function \eqref{eq:02.11} with $V(x) =\eta(\norm{x})$ corresponds to a Massera type of construction, which in its original form in \cite{Massera1949} is defined as:
\begin{equation}
\label{eq:02.13}
W(x(t))=\int_{t}^{\infty}\alpha(\norm{x(\tau)})\d \tau,
\end{equation}
with $\alpha\,:\,\R_{\geq 0}\,\rightarrow\, \R_{\geq 0}$ an appropriately chosen continuous function.
In \cite{Bjornsson2014} the construction above is recalled as
\begin{equation}
\label{eq:02.14}
W(x(0))=\int_{0}^{N}\norm{x(\tau)}_2 \d \tau,
\end{equation}
as an alternative to the construction in \cite[Theorem 4.14]{Khalil2002}, where the limits of the integral are $t$ and $t+N$.
This type of construction facilitated an extensive amount of converse results. Also, the converse proof in \cite{Massera1949} set up a proof technique which is based on the so--called Massera's Lemma \cite[p. 716]{Massera1949} for constructing the function $\alpha$ in \eqref{eq:02.13}.

The formulation in \cite{Bjornsson2014} is based on the exponential stability assumption, though an extension for the asymptotic stability  case is suggested. The extension relates to the construction in \eqref{eq:02.11}, where instead of the function $\eta$, a nonlinear scaling of the norm of the state trajectory is used. This scaling function is to obtained from a $\mc{K}\mc{L}$ estimate and Lemma~\ref{lemma:01.01} to provide an exponentially decreasing in time upper bound for the asymptotic stability estimate. In \cite{Bjornsson2015}, a similar construction to the one in \eqref{eq:02.11} is proposed, namely
\begin{equation}
\label{eq:02.14b}
W(x(0))=\int_{0}^{T}\gamma(x(\tau)) \d \tau, 
\end{equation}
with $T$ a positive, finite constant and $\gamma$ a positive definite function. While both the proposed construction in this paper and the one in \eqref{eq:02.14b} rely on arbitrary functions of the state, the main difference is in the choice of the integration interval, with $d$ being such that \eqref{eq:02.01b} holds for the arbitrary $\eta$ function. 

By using the finite--time function $V$, we provide an alternative to Massera's construction via Corollary~\ref{co:02.01},
whilst by defining $W$ as in \eqref{eq:02.11} we allow for a $\mc{K}\mc{L}$--stability assumption with the $\mc{K}\mc{L}$ function satisfying \eqref{eq:02.08}.
The freedom in choosing any function of the state norm $\eta\in\mc{K}_\infty$ in the proposed construction facilitates an implementable verification procedure which is detailed in Section~\ref{sec:04} and does not rely on a specific, possibly more complex form for a LF, which can add to the computational load.

 \section{Expansion scheme}
\label{sec:03}
In \cite{Chiang89} a scheme for constructing LFs starting from a given LF, which at every iterate provides a less conservative estimate of the DOA of a nonlinear system of the type \eqref{eq:01.01} was proposed. The sequence of Lyapunov functions of the type
\begin{eqnarray}
\label{eq:03.01}
\nonumber
W_1(x)&=&W(x+\alpha_1f(x))\\
\nonumber
W_2(x)&=&W_1(x+\alpha_2f(x))\\
&\vdots&\\
\nonumber
W_n(x) &=&W_{n-1}(x+\alpha_nf(x)),
\end{eqnarray}
with $\alpha_i\in\R_{\geq0}$, $i=1,2\ldots,n$ leads to the DOA estimates set inclusion
$$\mc{S}_W(c)\subset\mc{S}_{W_1}(c)\subset\ldots\subset\mc{S}_{W_n}(c),$$
where $\mc{S}_{W_i}(c):=\{x\in\R^n\,|\,W_i(x)\leq c\}$, $i=1,2\ldots,n$ denote the largest level sets of the Lyapunov functions generated by the expansion sequence \eqref{eq:03.01} with $\mc{S}_{W}(c)\subset\mc{S}$, where $\mc{S}$ is the $d$--invariant set.
Note that the largest level set of the LF $W$, included in the $d$--invariant set $\mc{S}$ is a subset of the true DOA of the system \eqref{eq:01.01}. Thus, the expansion scheme \eqref{eq:03.01} applied on a computed $W$ will provide a better estimate of the true shape of the DOA contained in $\mc{S}$. However, a less conservative initial set $\mc{S}$ leads to a better estimate of the true DOA.

We propose to utilize the expansion idea in \cite{Chiang89} to generate a sequence of FTLFs, with the purpose to generate a more appropriate $d$--invariant set. The next result follows as a consequence of Lemmas~$4$-$1$ and $4$-$2$ from \cite{Chiang 89}.
\begin{theorem}
\label{co:03.01}
Let $V$ be a FTLF, i.e. conditions \eqref{eq:02.01a} and \eqref{eq:02.01b} hold with respect to the $d$--invariant set $\mc{S}$. Furthermore let $V$ be continuously differentiable and define
$\mc{S}_V(c):=\{x\in\R^n\,|\,V(x)\leq c\}\subset\mc{S}$.
Then, there exists an $\alpha\in\R_{\geq0}$ such that for $V_1(x)=V(x+\alpha f(x))$ and $\mc{S}_{V_1}(c):=\{x\in\mb{R}^n\,|\,V_1(x)\leq c\}$, it holds that
$\mc{S}_V(c)\subset\mc{S}_{V_1}(c)$ and $V(x+\alpha f(x))$ is a FTLF.
\end{theorem}
\begin{IEEEproof}
By hypothesis, $V(x)$ is a continuously differentiable FTLF, i.e.  let $V$ have continuous partial derivatives of order $r$ higher or equal to $1$.  Then by \cite[p.190]{Dieudonne69} and \cite[Lemma 2-1]{Chiang89} we can write:
\begin{equation}
\nonumber
\begin{split}
V_1(x(t)) =&\,V(x(t))+\alpha\dot{V}(x(t)) + \frac{\alpha^2}{2!}\ddot{V}(x(t))+\ldots+\\
&\,\frac{\alpha^{r-1}}{(r-1)!}V^{(r-1)}(x(t)) +\frac{\alpha^{r}}{r!}V^{(r)}(z(t)),
\end{split}
\end{equation}
where $V^{r}$ denotes the $r$-th derivative of $V$, i.e. $V^{(r)}(x)=\nabla^{(r-1)} V(x)^\top f(x)$ and $z(t)=x(t)+h\alpha f(x(t))$, for $h\in[0,1]$.

Let $V(x(t))= (V\circ x)(t)=\psi(t)$.
Similarly, by the formula \cite[p.190, (8.14.3)]{Dieudonne69}, we have that
\begin{equation}
\nonumber
\begin{split}
\psi(t+d) = &\,\psi(t)+d\dot{\psi}(t)+\frac{d^2}{2!}\ddot{\psi}(t) +\dots
\frac{d^{r-1}}{(r-1)!}\psi^{(r-1)}(t) + \\
&\,\frac{d^{r}}{r!}\psi^{r}(w),
\end{split}
\end{equation}
where $w = t+hd$, $h\in[0,1]$ and $\dot{\psi}(t)=\dot{V}(x(t))=\nabla V(x)^\top f(x(t))$.
Then,
\begin{equation}
\nonumber
\begin{split}
V(x(t+d))=&\,V(x(t))+d\dot{V}(x(t)) +\frac{d^2}{2!}\ddot{V}(x(t)) +\ldots +\\ 
&\,\frac{\alpha^{r-1}}{(r-1)!}V^{(r-1)}(x(t)) + \frac{\alpha^{r}}{r!}V^r(x(w)).   
\end{split}
\end{equation} 
From the expressions of $V_1(x(t))$ and $V(x(t+d))$, we can write
\begin{equation}
\label{eq:03.02}
\begin{split}
V_1(x(t)) \leq &\,V(x(t))+V(x(t+d))-V(x(t))+ \\
&\,|\alpha-d|\dot{V}(x(t)) + \frac{|\alpha^2-d^2|}{2!}\ddot{V}(x(t)) + \ldots +\\
&\,  \frac{|\alpha^{r-1}-d^{r-1}|}{(r-1)!}V^{(r-1)}(x(t)) +\\
&\,\bigg|\frac{\alpha^r}{r!}V^{r}(z(t)) -
 \frac{\alpha^r}{r!} V^{r}(x(w))\bigg|\\
\leq &\, V(x(t))+V(x(t+d))-V(x(t))+\\
&\, |\alpha-d|\dot{V}(x(t)) +|\alpha-d|\epsilon.
\end{split}
\end{equation}

Since $V(x(t+d))-V(x(t))<0$, there exists a $\beta\in\R_{> 0}$ such that $V(x(t+d))-V(x(t))<-\beta$, for all $x(0)\in\mc{S}\setminus\mc{S}_V(c)$, where $\mc{S}\setminus\mc{S}_V(c)$ is a compact set containing no equilibrium point. Furthermore, on the compact set  $\mc{S}\setminus\mc{S}_V(c)$, $\epsilon>0$ is a bound on the sum of higher order continuous terms in the above expression. This is due to the fact that continuous functions are bounded on compact sets. As such, we obtain that
\begin{equation}
\nonumber
\begin{split}
V_1(x(t)) < & \,V(x(t)) -\beta +|\alpha- d|
         (\dot{V}(x(t))+ \epsilon)\\
\end{split}
\end{equation}
Let $|\alpha-d|<\bar{\alpha}$  and let $\dot{V}(x(t))+\epsilon\leq\nu,$ $\nu\in\R_{>0}$ for all $x\in\mc{S}\setminus\mc{S}_V(c)$ since every continuous function is bounded on a compact set. Thus, we obtain that
$$V_1(x(t))<V(x(t)) +\bar{\alpha}\nu-\beta.$$
If $(\bar{\alpha}\nu-\beta)<0$, hence for $0<|\alpha-d|<\beta/\nu$, it holds that
$$V_1(x(t))<V(x(t)).$$
Similarly as in \cite[Lemma $4$-$1$]{Chiang89}, this implies that $\mc{S}_V(c)\subset\mc{S}_{V_1}(c)$.

Next we will show that $V_1$ is a FTLF. 
From Lemma~\ref{lemma:02.01} we know that the function  $W(x(t))=\int_t^{t+d} V(x(\tau))\d\tau$ is a LF for \eqref{eq:01.01}. From \cite[Lemma 4-2]{Chiang89} it is known that there exists some $\alpha>0$ such that 
$W_1(x(t))= W(x(t)+\alpha f(x(t)))=\int_t^{t+d} V(x(\tau)+\alpha f(x(\tau)))\d\tau$
is a LF. This implies that there exists some $\mc{K}$--function $\tilde{\gamma}(x(t))$ such that
\begin{equation}\nonumber
\begin{split}
\dot{W}_1(x(t))&\,=V(x(t+d) +\alpha f(x(t+d)))-\\
&\quad \,\,\,   V(x(t)+\alpha f(x(t)))\\
&\,\leq-\tilde{\gamma}(x(t)),
\end{split}
\end{equation}
where the Leibniz integral rule was used. This further implies that
$V_1(x(t+d))-V_1(x(t))\leq-\tilde{\gamma}(x(t))$, thus $V_1$ is a FTLF. 
\end{IEEEproof}

\section{Verification}
\label{sec:04}
\subsection{The indirect approach}
\label{sec:04.01}
In Section~\ref{sec:02} it has been shown that if the equilibrium of a given system is $\mc{K}\mc{L}$--stable, then a method to construct a Lyapunov function is provided by \eqref{eq:02.11}, for $V(x)$ defined by any function $\eta\in\mc{K}_\infty$ and any norm.
The method is constructive starting with a given candidate  $d$--invariant set $\mc{S}$ and a candidate function $V(x)=\eta(\norm{x})$. Due to the $d$--invariance property of $\mc{S}$, verifying condition \eqref{eq:02.01b} for the chosen $V$ is reduced to verifying
\begin{equation}
\label{eq:04.05}
V(x(d))-V(x(0))\leq-\gamma(\norm{x})<0,
\end{equation}
for all $x(0)\in\mc{S}$. The difficulty in verifying \eqref{eq:04.05} is given by the need to compute $x(d)$, for all $x(0)\in\mc{S}$. However if $x(d)$ is known analytically, then it suffices to verify \eqref{eq:04.05}  for all initial conditions in a chosen set $\mc{S}$.  Then the largest level set of $W$, defined as in \eqref{eq:02.11}, which is included in $\mc{S}$, is a subset of the true DOA of the considered system. The verification of \eqref{eq:04.05} translates into solving a problem of the type
\begin{equation}
\label{eq:04.05.ver}
\begin{split}
&\max_{x(0)}  \,[V(x(d))-V(x(0))]\\
&\mbox{subject to } x(0)\in\mc{S}.
\end{split}
\end{equation}
The  regularity assumptions on the map describing the dynamics \eqref{eq:01.01} ensures that the solution is continuous for all $t\in[0,d]$. Furthermore, since $V(x)=\eta(\norm{x})\in\mc{K}_\infty$ and the $d$--invariant candidate set $\mc{S}$ is compact, the problem \eqref{eq:04.05.ver} will always have a global optimum. A way to avoid solving \eqref{eq:04.05.ver} will be indicated in Section~\ref{sec:04.02}.

When the analytical solution is not known, or obtaining a numerical approximation is computationally tedious, as it can be the case for higher order nonlinear systems, then we propose the following approach starting from the linearized dynamics of \eqref{eq:01.01}. In what follows we recall some relevant properties of the map $f$ with respect to its linearization. The detailed derivations can be found in \cite[Chaper 4.3]{Khalil2002}.
Firstly, by the mean value theorem it follows that
\begin{equation}
\label{eq:04.01}
f_i(x)=f_i(0)+\frac{\partial f_i}{\partial x}(z_i) x,\quad i =1,\ldots,n,
\end{equation}
where $z_i$ is a point on the line connecting $x$ to the origin. Since the origin is an equilibrium of \eqref{eq:01.01},
$$f_i(x)=\frac{\partial f_i}{\partial x}(z_i) x = \frac{\partial f_i}{\partial x}(0)x + \left(\frac{\partial f_i}{\partial x}(z_i)-\frac{\partial f_i}{\partial x}(0)\right)x.$$
Then,
\begin{equation}
\label{eq:04.02}
f(x)=A x +g(x),
\end{equation}
with
$$A = \frac{\partial f}{\partial x}(0)=\left[\frac{\partial f(x)}{\partial x}\right]_{x=0},$$ $$g_i(x)=\left(\frac{\partial f_i}{\partial x}(z_i)-\frac{\partial f_i}{\partial x}(0)\right)x.$$
Furthermore,
$$\norm{g_i(x)}\leq \norm{\frac{\partial f_i}{\partial x}(z_i)-\frac{\partial f_i}{\partial x}(0)}\norm{x}$$ and
$$\frac{\norm{g(x)}}{\norm{x}}\rightarrow 0,\, \mbox{as}\, \norm{x}\rightarrow 0.$$
In this way, in a sufficiently small region around the origin we can approximate the system \eqref{eq:01.01} with $\dot{x}=Ax$.

The next result is an analog to \emph{Lyapunov's indirect method} \cite[Theorem 4.7]{Khalil2002}, but in terms of the FTLF concept. The aim is to provide a validity result for the FT condition \eqref{eq:02.01b} for the nonlinear system \eqref{eq:01.01}, whenever a (global) FT type condition is satisfied for the linearized system with respect to the origin.
\begin{theorem}
\label{th:04.01}
Let $V(x)=\norm{x}$ be a global FTLF function for  $\dot{x}= Ax$, i.e. there exists a  $d>0$ such that $\norm{e^{Ad}} < 1$. Additionally, let
\begin{equation}
\label{eq:04.03}
e^{d\mu(A)}-1=-\varsigma,\quad \varsigma\in\R_{>0}
\end{equation}
be satisfied.
Then the following statements hold.
\begin{enumerate}[1.]
\item  There exists a $d$--invariant set $\mc{S}$ for which $V(x)$ is a FTLF for \eqref{eq:01.01}.
\item  There exists a set $\mc{A}\subseteq\mc{S}$ for which
\begin{equation}
\label{eq:04.W}
W(x) = \int_0^d V(x+\tau f(x))\d\tau
\end{equation}
is a LF for \eqref{eq:01.01}.
\end{enumerate}
\end{theorem}
\begin{IEEEproof}
We start by proving point 1.
As indicated in \cite{Khalil2002}, from $\frac{\norm{g(x)}}{\norm{x}}\rightarrow 0$, as $\norm{x}\rightarrow 0$, it follows that for any $\delta>0$, there exists an $r>0$ such that $\norm{g(x)}<\delta\norm{x}$, for $\norm{x}<r$. Thus, the solution of the system defined with the map in \eqref{eq:04.02} is bounded, whenever $g(x)$ is bounded \cite{Soderlind}, as shown below.
\begin{equation}
\nonumber
\begin{split}
\norm{x(d)}\leq &\, e^{d\mu(A)}\norm{x(0)} +\int_0^d e^{(d-\tau)\mu(A)}\norm{g(x(\tau))}\d\tau\\
           \leq &\, e^{d\mu(A)}\norm{x(0)} + \int_0^d e^{(d-\tau)\mu(A)}\delta\norm{x(\tau)}\d\tau.
\end{split}
\end{equation}
By applying the Bellman--Gronwall Lemma~\ref{lemma:01.03} to the inequality above, we obtain
\begin{equation}
\nonumber
\begin{split}
\norm{x(d)} \leq &\, e^{d\mu(A)}\norm{x(0)}e^{\int_0^d e^{(d-\tau)\mu(A)}\delta\d\tau}\\
            \leq &\, e^{d\mu(A)}\norm{x(0)}e^{\delta\frac{e^{d\mu(A)}-1}{\mu(A)}}\\
            \leq &\,e^{d\mu(A)-\frac{\varsigma\delta}{\mu(A)}}\norm{x(0)}.
\end{split}
\end{equation}
Thus, $V(x(d))\leq \rho V(x(0))$, with $\rho := e^{d\mu(A)-\frac{\varsigma\delta}{\mu(A)}}$. For the equivalent FT condition \eqref{eq:02.06} to hold, $\rho$ must be subunitary, and equivalently
\begin{equation}
\label{eq:04.04}
d\mu(A)-\frac{\varsigma\delta}{\mu(A)}<0.
\end{equation}
For equation \eqref{eq:04.03} to hold, $\mu(A)$ must be negative. Thus, we obtain that
\begin{equation}
\label{eq:04.bound}
\delta\leq\frac{d\mu(A)^2}{\varsigma},
\end{equation}
which provides an upper bound on $\delta$, and thus on $r$. Consequently, there exists a $d$--invariant set $\mc{S}\subseteq\{x\in\R^n\,|\,\norm{x}< r\}$.

As for the second item of the theorem, let us consider the Dini derivative expression for $\dot{W}(x(0))$,
$$\dot{W}(x(0))=\diniD^+ W(x(0))= \limsup_{h\rightarrow 0^+}\frac{W(x(h))-W(x(0))}{h}.$$
Then,
\begin{equation}
\nonumber
\begin{split}
\diniD^+ W(x(0))
= &\,\limsup_{h\rightarrow 0^+}\bigg[\frac{\int_h^{h+d} V(x(0)+\tau f(x(0)))\d\tau-}{h} \\
&\frac{\int_0^d V(x(0)+\tau f(x(0)))\d\tau}{h}\bigg]\\
=&\,\limsup_{h\rightarrow 0^+}\bigg[\frac{\int_d^{d+h} V(x(0)+\tau f(x(0)))\d\tau - }{h}\\
&\frac{\int_0^{h} V(x(0)+\tau f(x(0)))\d\tau}{h}\bigg]\\
=\footnotemark &\, V(x(0)+df(x(0)))-V(x(0))\\
< &\, V(x(d))-V(x(0))\leq-\gamma(\norm{x(0)}),
\end{split}
\end{equation}
\footnotetext{Here we applied L'Hospital rule together with Leibniz integral rule. }
where we used the fact that $V_1(x(t))<V(x(t))$  shown in the proof of Corollary~\ref{co:03.01}.
\end{IEEEproof}

In the theorem above,  due to the equivalence result in Lemma~\ref{lemma:02.01}, existence of a FTLF $V$ is equivalent to existence of a true LF $W$ defined as in \eqref{eq:02.11}. Since $V$ is only valid in the region around the origin defined by $\delta$, $W$ will also be a valid LF in some subset of that region. However, the expression of $W$ in \eqref{eq:02.11} still involves knowing the solution $x(d)$. In this case, relying on the solution of the linearized system might lead to conservative approximations of the DOA. In view of the fact that we want to construct LFs that lead to relevant DOA estimates, the construction in \eqref{eq:04.W} is more suitable as it includes the  nonlinear vector field. 

The condition in \eqref{eq:04.03}, essentially requires that there should exist a (weighted) norm for which the induced logarithmic norm of the matrix obtained by evaluating the Jacobian associated to \eqref{eq:01.01} at the origin is negative.   The choice of the norm inducing the matrix measure is dictated by the choice of the norm defining the FTLF $V(x)$.  A similar condition has been introduced in \cite{Sontag2014} for characterizing infinitesimally contracting systems on a given convex set $\mb{X}\subseteq\R^n$. Therein, the logarithmic norm of the Jacobian at all points in the set $\mb{X}$ is required to be negative.
\subsection{Computational procedure}
\label{sec:04.02}
\subsubsection{Compute a $d$ for which the finite--time condition holds for the linearized system with respect to the origin}
Let
$$\dot{\delta} x = \left[\frac{\partial f(x)}{\partial x}\right]_{x=0}\delta x,$$
be the linearized system. Then condition \eqref{eq:04.05} can be verified as:
$$\eta(\norm{e^{d\left[\frac{\partial f(x)}{\partial x}\right]_{x=0}}x(0)})-\eta(\norm{x(0)})<0,$$
for all $x(0)$ in some compact, proper set $\mc{S}$. As such, for a given value of $d$, \eqref{eq:04.05} can be verified via a feasibility problem. If $\eta =\id$ then the feasibility problem above translates into verifying the matrix norm condition:
\begin{equation}
\label{eq:norm}
 \norm{e^{d\left[\frac{\partial f(x)}{\partial x}\right]_{x=0}}}<1.
\end{equation}
Note that if there exists a $d>0$ such that the condition \eqref{eq:04.03} holds for some $\varsigma>0$, then condition \eqref{eq:norm} is implicitly satisfied since it holds that $\norm{e^{dA}}\leq e^{d\mu(A)}$ for some real matrix $A$ \cite{Soderlind}.
\subsubsection{ Compute $W(x)$}
\begin{equation}
\label{eq:04.06}
W(x)=\int_{0}^d V(x+\tau f(x))\d\tau.
\end{equation}
In this construction we exploit the results in Corollary~\ref{co:03.01} and Theorem~\ref{th:04.01}, which guarantee that $V(x+\tau f(x))$ remains a FTLF for \eqref{eq:01.01} for all $\tau\in[0,d]$.
\subsubsection{Find the best DOA estimate of \eqref{eq:01.01} provided by $W$}
Let $C$ be such that the set $\{x\in\R^n\,|\,W(x)\leq C\}$ is included in the set
$\{x\in\R^n\,|\,\nabla^\top W f(x)=0\}$. Finding $C$, implies solving an optimization problem, which involves rather complex nonlinear functions. However feasibility problems (for example by using bisection) can be solved successfully to obtain the best $C$ which leads to a true DOA estimate.
The feasibility problem is as follows:
\begin{equation}
\label{eq:04.07}
\begin{split}
&\max_{x}  \,\nabla^\top W(x) f(x)\\
&\mbox{subject to } W(x)\leq C,
\end{split}
\end{equation}
for a given $C$ value. The largest $C$ value for which $\dot{W}(x)$ remains negative renders the best DOA estimate provided by a true LF $W(x)$, which is valid for \eqref{eq:01.01} since $\dot{W}(x)=\nabla^\top W(x) f(x)$. The optimization problem \eqref{eq:04.07} is standard in checking validity regions of LFs aposteriori to construction. For more details, see \cite{Chesi2011} or \cite{Hachicho2007}.

The bound \eqref{eq:04.bound}, via $r$,  provides an implicit, a priori theoretical indication of the region where $W$ is a valid Lyapunov function, or in other words of a subset of the DOA. An explicit estimate can be obtained a posteriori to computing $W$ by solving the problem above.
\subsubsection{Further improve the DOA estimate by expansion of $W$}
The DOA estimate obtained above can be further improved by making use of expansion methods as introduced in \cite{Chiang89}. More specifically, the expansion method states that the $C$--level set of $W_1(x)=W_1(x+\alpha f(x))$ for some $\alpha\in\R_{>0}$ is a subset of the true DOA of the considered system and it will contain the estimate of the DOA obtained by the $C$--level set of $W$.

In summary, the proposed method starts by verifying the finite--time decrease condition \eqref{eq:04.05} for a candidate function $V(x)=\eta(\norm{x})$, $\eta\in\mc{K}_\infty$ and a candidate $d$--invariant set $\mc{S}$. The simplest way to do this, while avoiding solution approximations,  given in the first step of the computation procedure above, is to verify \eqref{eq:04.05} globally ($x\in\R^n$) for the linearization of \eqref{eq:01.01} around the origin.  Then, it is known by Theorem~\ref{th:04.01}, that there exists a set $\mc{S}\subseteq\{x\in\R^n\,|\,\norm{x}<r\}$ which is $d$--invariant, such that \eqref{eq:04.05} holds for the true nonlinear system.

Next, in the second step of the procedure,  $W$ is computed via the analytic formula \eqref{eq:04.06}, which yields an educated guess of a true LF. Thus, the final check in the third step is a verification of the standard Lyapunov condition on $W$, with $W$ known, and with the aim to maximize the level set $C$ of $W$ where the condition is satisfied.

In contrast,  most of the other proposed methods in the literature, compute the LF $W$ simultaneously with verifying the derivative negative definiteness condition, which is in a general a more difficult task, even because it is not clear how to select a non--conservative $W$.

\section{Examples}
\label{sec:05}
\subsection{An illustrative example with no polynomial LF}
\label{ex:05.01}
Consider the system
\begin{eqnarray}
\label{eq:05.01}
\dot{x}_1&=&-x_1+x_1x_2\\
\nonumber
\dot{x}_2&=&-x_2,
\end{eqnarray}
with solutions
$x_1(t)=x_1(0)e^{(x_2(0)-x_2(0)e^{-t} -t)}$ and  $x_2(t)=x_2(0)e^{-t}$.
In \cite{Parillo2011} it has been shown that the system is GAS by using the Lyapunov function
$V_{GAS}(x)=\ln{(1+x_1^2)} +x_2^2$. Furthermore, it has been shown that no polynomial LF exists for this system. We will illustrate our converse results on this system and we will show that the proposed constructive approach leads to local approximations of the DOA that cover a larger area of the state space compared to similar sublevel sets of $V_{GAS}$.
Since the system is GAS, then it is $\mc{K}\mc{L}$--stable with $\mc{S}=\R^n$, and  by Theorem~\ref{th:02.02} we have that for any function $\eta\in\mc{K}_\infty$, there exists a set $\mc{S}\subseteq\R^n$ and a scalar $d$ such that the function $V(x)=\eta(\norm{x})$, for any $x\in\mc{S}$ is a FTLF for the system \eqref{eq:05.01}.
\begin{figure}[htp]
\begin{center}
\subfigure[The set $\mb{X}$--black, the set $\mc{S}$ defined by the level set  $C_V=1.6$ of FTLF $V$--red and the level set defined by $C_W=0.415$ of the LF $W$--blue together with the vector field plot of \eqref{eq:01.01}.]
{\includegraphics[width=5.1cm,height=4.7cm]{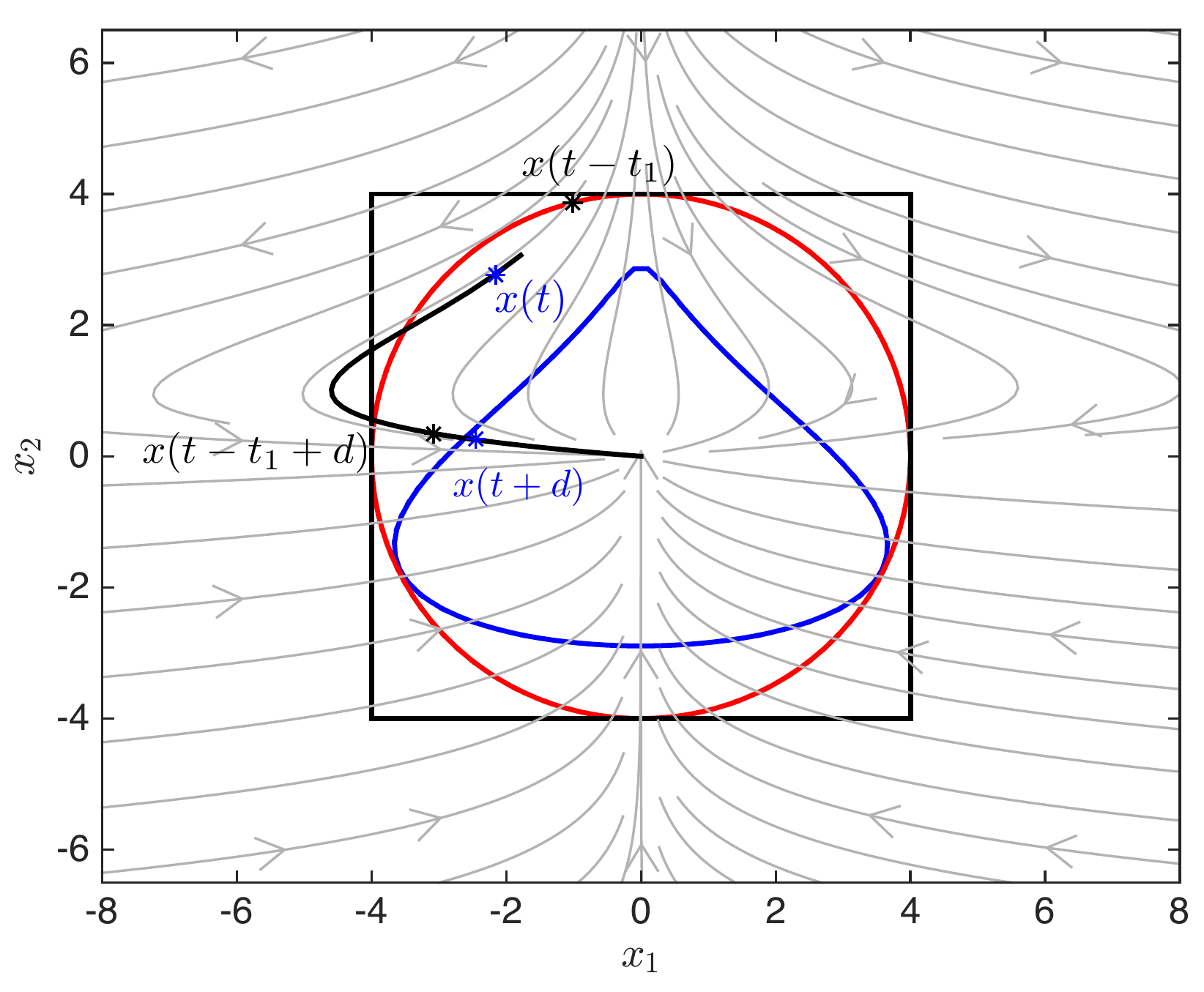}
\label{fig:05.01a}}
\hspace{-0.0cm}
\subfigure[The same sets as in Figure~\ref{fig:05.01a} together with the level set defined by $C_{W_1}=1.5$ of  $W_1$--green.]{\includegraphics[width=5.1cm,height=4.55cm]{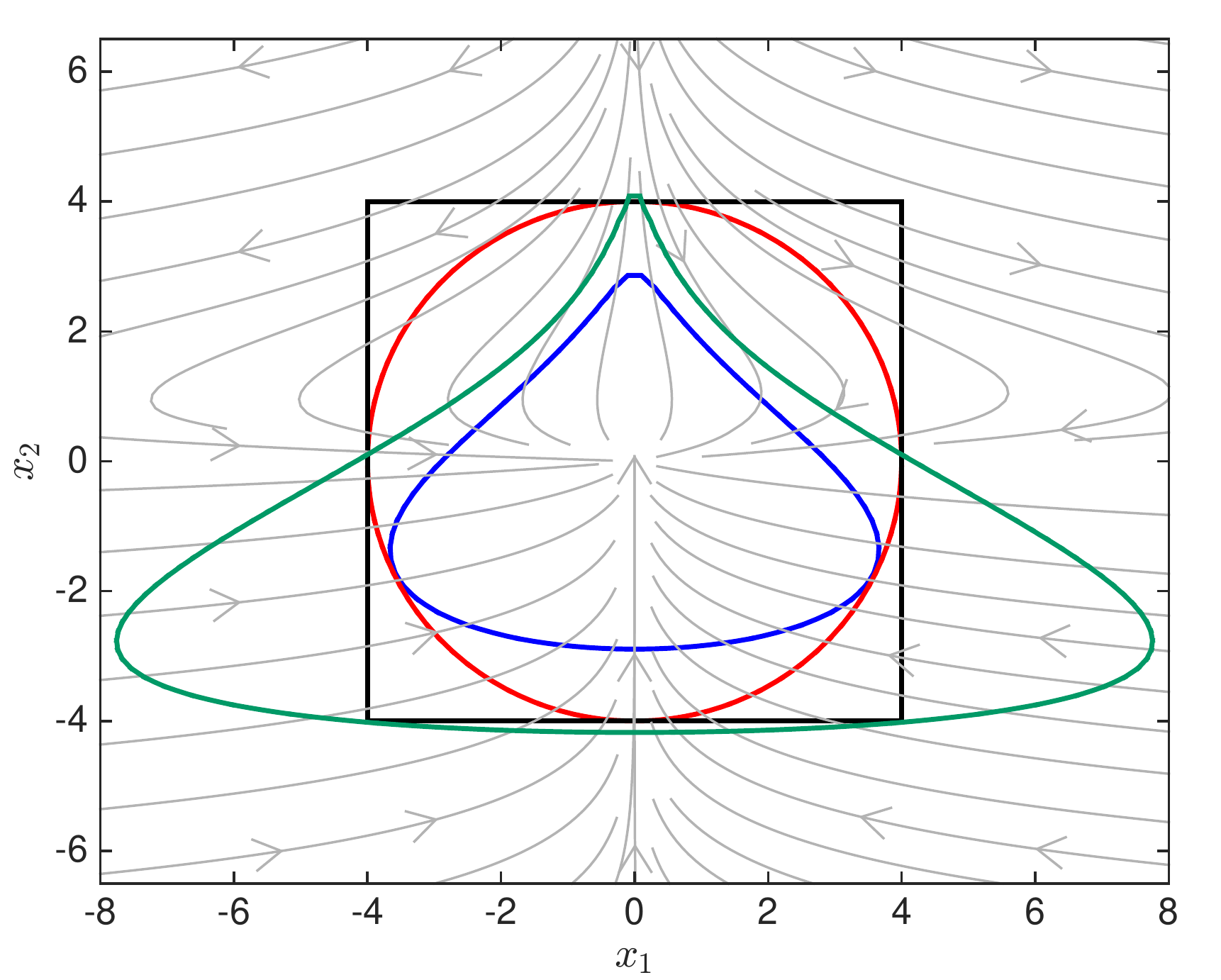}
\label{fig:05.01b}}
\subfigure[Comparison plot: $W(x)=0.415$--blue, $W_1(x)=1.5$--green and $V_{GAS}=3$--gray.]
{
\includegraphics[width=5.1cm,height=4.7cm]{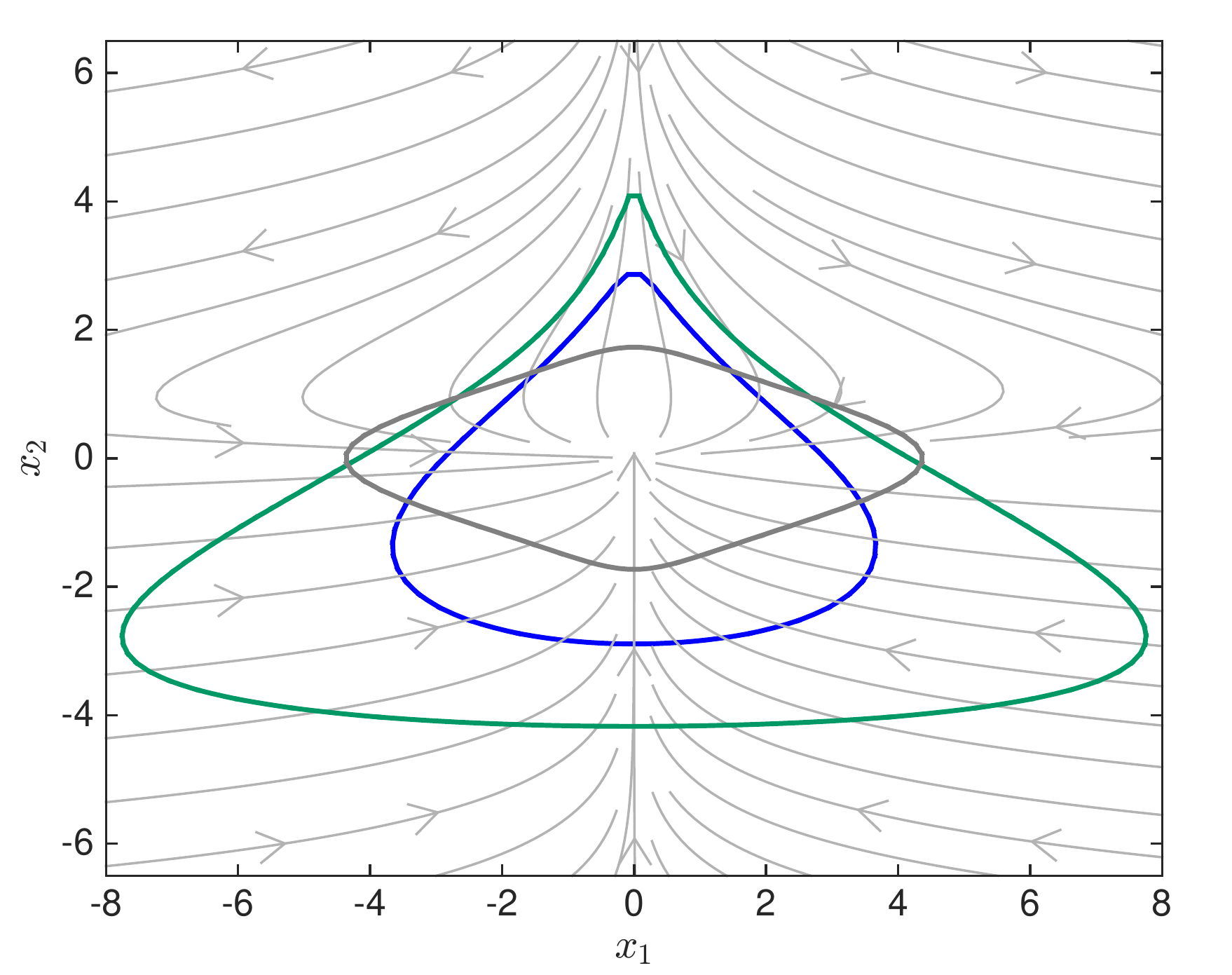}
\label{fig:05.01c}}
%
\caption{Plots of the level sets of the function $V$, $W$ and $W_1$ computed for \eqref{eq:01.01}.}
\label{fig:05.01}
\end{center}
\end{figure}
Consider the set $\mb{X}:=\{x\in\R^n\,|\,\norm{x}_\infty\leq 4\}$, which is displayed in Figure~\ref{fig:05.01} with the black contour. Pick $V(x)=x^\top P x$,
where $P=\left(\begin{smallmatrix}   0.1 & 0\\
                                    0 & 0.1
                                   \end{smallmatrix}\right).$
For this choice of $V$, the feasibility problem \eqref{eq:04.05.ver} was solved using the \verb|sqp| algorithm with \verb|fmincon|. Since all involved functions are Lipschitz and we considered a polytopic, compact candidate set $\mc{S}$,  the problem \eqref{eq:04.05.ver} has a global optimum.
Thus condition \eqref{eq:04.05} and consequently, \eqref{eq:02.01b} holds for $d=2.4$ for any $x\in\mc{S}$ where $\mc{S}$ is the largest level set of $V$ included in the set $\mb{X}$; $\mc{S}$ is shown in Figure~\ref{fig:05.01a} in red.
By Lemma~\ref{lemma:02.01}, we obtain that
$$W(x(t))=\int_{t}^{t+d}x(\tau)^\top P x(\tau) \d\tau$$ is a Lyapunov function for \eqref{eq:05.01},
for any $t\in\R_{\geq 0}$.
From \eqref{eq:04.05} and the equivalence result in Lemma~\ref{lemma:02.01}, it is sufficient to compute $W(x(0))$.
Since the system is GAS, any scaling of the set $\mc{S}$ will satisfy \eqref{eq:02.01b}, however with a bigger $d$, due to the nonlinear dynamics. For the computed $d$ and chosen $\mc{S}$ and $V$, let $\mc{S}_V(C_V)$  denote the largest level set of $V$ included in $\mc{S}$.
Then, the largest level set of  $W$ in $\mc{S}_V(C_V)$ will be a subset of the true DOA of the system. In Figure~\ref{fig:05.01a} we show the level set of $W$ defined by $C_W=0.415$ in blue.
Next, we will illustrate also the expansion method.
Consider $W_1(x(t))=W(x(t)+\alpha_1f(x(t))$, with $\alpha_1=0.1$. The level set defined by $W_1(x) =C_{W_1}$, where $C_{W_1} = 1.5$ is shown in Figure~\ref{fig:05.01b} in green. Note that, by construction, $W_1(x)=1.5$ is not restricted to the black set anymore. For the sake of comparison we show in Figure~\ref{fig:05.01c} a plot of level sets of the two computed functions, $W$ and $W_1$ and a relevant level set value of the logarithmic LF $V_{GAS}$.

\subsection{3D example from literature}
\label{ex:05.02}
We consider the system in \cite[Example 5]{SiggiMicnon2015}, described by
\begin{eqnarray}
\label{eq:05.02}
\nonumber
\dot{x}_1&=&x_1(x_1^2+x_2^2-1)-x_2(x_3^2+1)\\
\dot{x}_2&=&x_2(x_1^2+x_2^2-1)+x_1(x_3^2+1)\\
\nonumber
\dot{x}_3&=&10x_3(x_3^2 -1).
\end{eqnarray}
\begin{figure}[htp]
\centerline{\includegraphics[width=0.45\columnwidth]{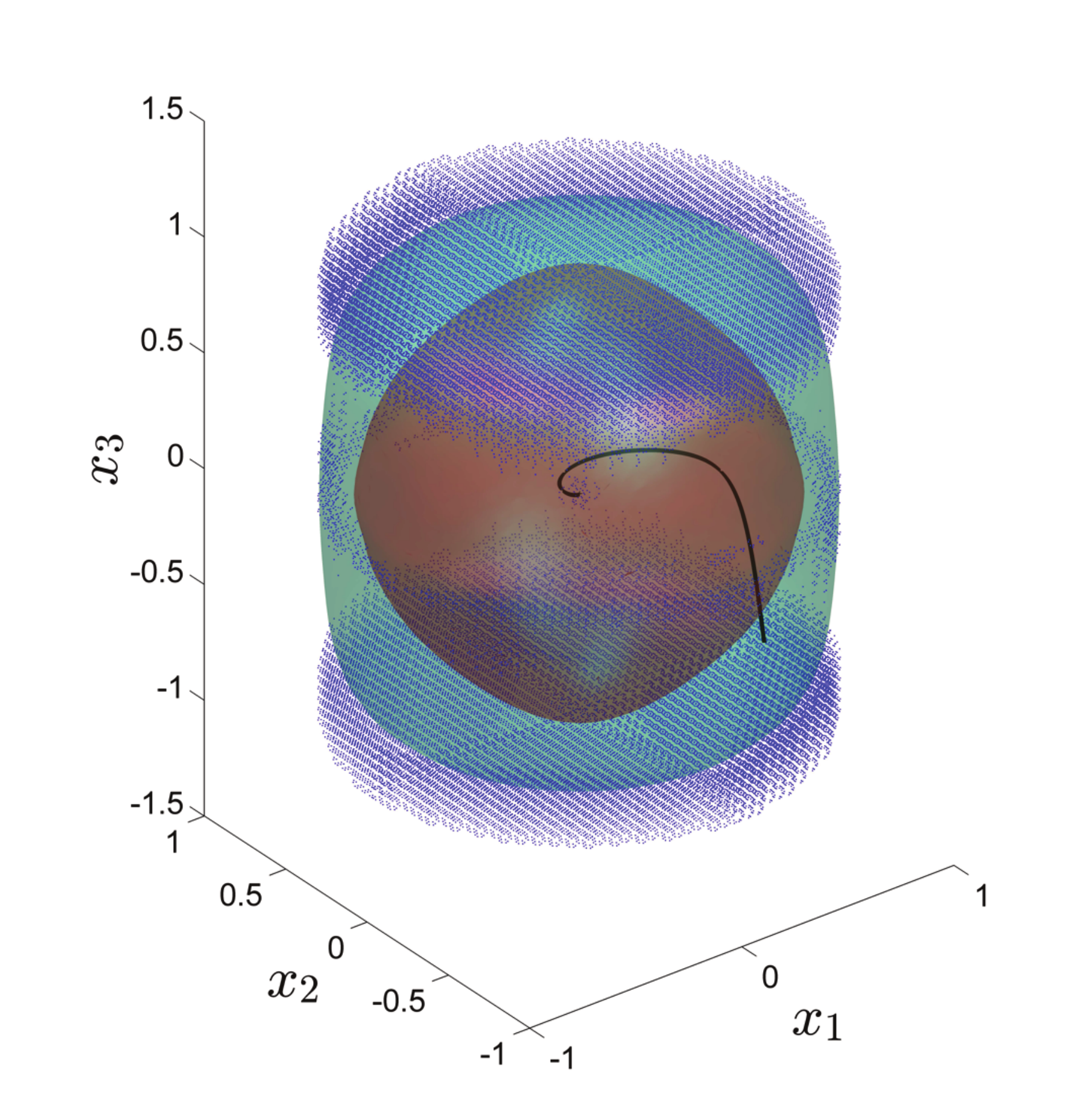}}
\caption{Plot of the the level set $W(x)=0.19$--green and the DOA approximation from \cite{SiggiMicnon2015}--red.}
\label{fig:05.02a}
\end{figure}
%
\begin{figure}[htp]
\centerline{\includegraphics[width=5.0cm,height=3.9cm]{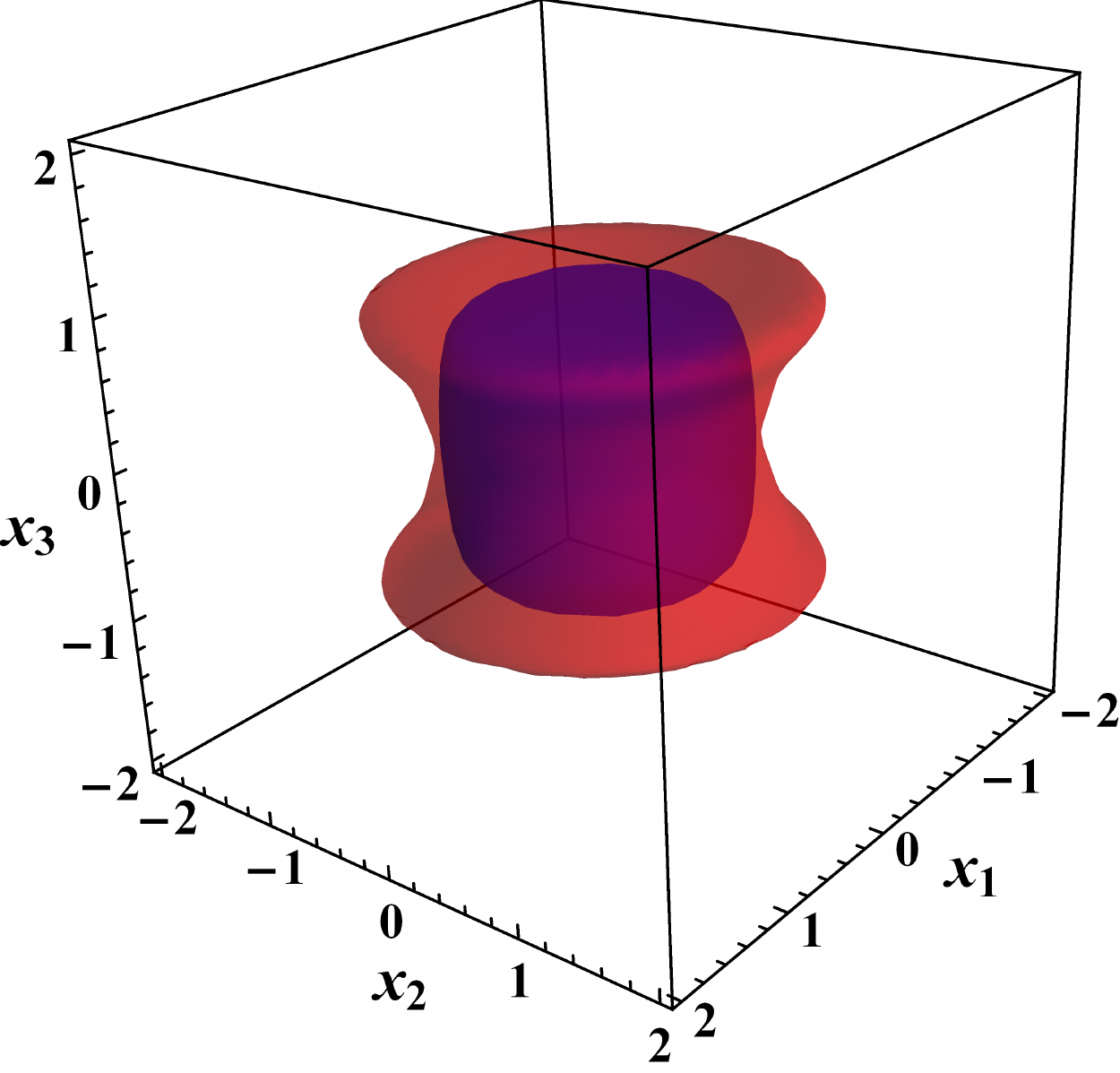}}
\caption{Validation of level set of $W$: $\nabla^\top W f=0$--red, level set of $W(x)=0.19$--green.}
\label{fig:05.02b}
\end{figure}
\vspace{-0.0cm}
For this system a piecewise affine LF was computed in \cite{SiggiMicnon2015}, resulting in the DOA plotted in Figure~\ref{fig:05.02a} with red. This plot can be found at \emph{www.ru.is/kennarar/sigurdurh/MICNON2015CPP.rar}. The blue dots in the figure represent the infeasibility points from computing the piecewise affine LF.

For this system, we do not know an analytic expression for $x(d)$. Thus, we apply the steps described in Section~\ref{sec:04.02}  for verification starting with the linearized dynamics. Let $V(x)=x^\top P x$ where $P$ is the identity matrix. Then, the condition \eqref{eq:norm} holds with $d=0.2$.
In Figure~\ref{fig:05.02b} the level set defined by $W(x)=0.19$ is plotted with green (inner set) together with the zero level set of $\dot{W}(x)$ in red. The value of $C=0.19$ was obtained by solving the feasibility problem \eqref{eq:04.07}. In Figure~\ref{fig:05.02a} we show the level set $W(x)=0.19$ and a trajectory of the system \eqref{eq:05.02}, initialized in $x_0=[0.5933\,-0.3636\,-0.6869]^\top$ with $W(x_0)=0.1198$, which is one of the infeasible points in \cite{SiggiMicnon2015}.
As for computation time, the example was worked out using MatlabR2015b, on a MacBookPro 2,8GHz Intel Core i5 processor and resulted in total computation time of $21.2712$s, while solving the problem \eqref{eq:04.07} takes $15.2496$s using \verb|fmincon| with the optimization algorithm \verb|sqp|. Finding the expression of $W$ takes $0.1726$s.

\subsection{Nonpolynomial 2D system--the genetic toggle switch}
\label{ex:05.03}
Consider the genetic toggle switch in Escherichia coli constructed in \cite{Gardner2000},
\begin{eqnarray}
\label{eq:05.03}
\nonumber
\dot{x}_1&=&\frac{\alpha_1}{1+x_2^\beta}-x_1\\
\dot{x}_2&=&\frac{\alpha_2}{1+x_1^\gamma}-x_2.
\end{eqnarray}
\begin{figure}[htp]
\centerline{\includegraphics[width=6.5cm,height=4.8cm]{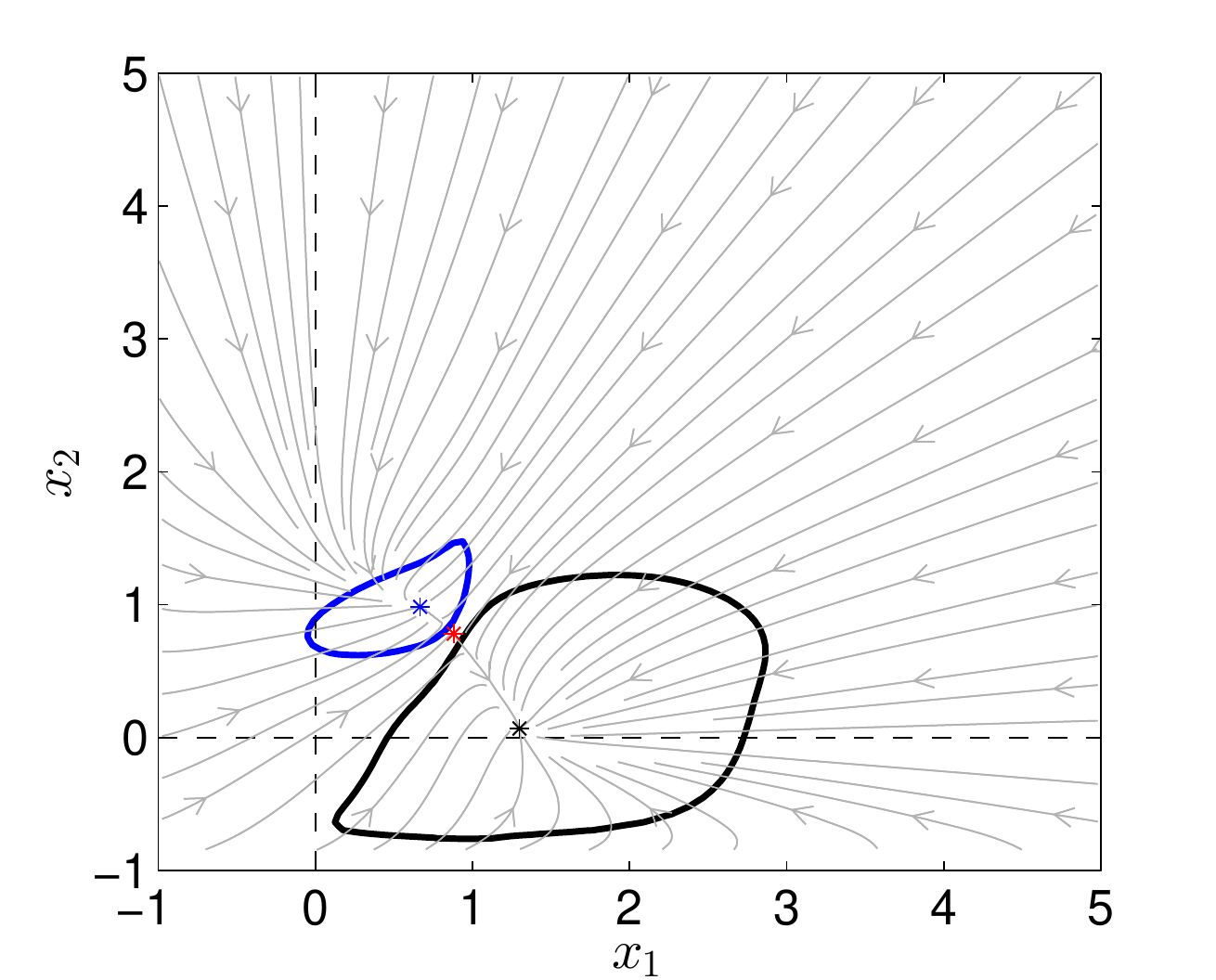}}
\caption{Level sets of the computed LFs corresponding to $E_1$--blue and $E_3$--red for the toggle switch system \eqref{eq:05.03} together with vector field plots.}
\label{fig:05.03}
\end{figure}
The genetic toggle switch is a synthetic, bistable gene--regulatory network, constructed from any two repressible promoters arranged in a mutually inhibitory network. The model \eqref{eq:05.03}, proposed in \cite{Gardner2000}, was derived from a biochemical rate equation formulation of gene expression. For this model, the set of parameters for which bistability is ensured is especially of interest, as it accommodates the real behavior of the toggle switch. The implications of the toggle switch circuit as an addressable memory unit are in biotechnology and gene therapy.
When at least one of the parameters $\beta,\gamma>1$, bistability occurs. In \eqref{eq:05.03}, $x_1$ denotes the concentration of repressor $1$, $x_2$ is the concentration of repressor $2$, $\alpha_1$ is the effective rate of synthesis of repressor $1$, $\alpha_2$ is the effective rate of synthesis of repressor $2$, $\beta$ is the cooperativity of repression of promoter $2$ and $\gamma$ is the cooperativity of repression of promoter $1$. The rational terms in the above equations represent the cooperative repression of constitutively transcribed promoters and the linear terms represent the degradation/dilution of the repressors.

The two possible stable states are one in which promoter $1$ transcribes repressor $2$ and one in which promoter $2$ transcribes repressor $2$. Let the parameters be defined by $\alpha_1=1.3$, $\alpha_2=1$, $\beta=3$ and $\gamma=10$. For this set of parameters the stable equilibria are
$E_1=\begin{pmatrix} 0.668 & 0.9829 \end{pmatrix}$ and $E_3=\begin{pmatrix} 1.2996 & 0.0678 \end{pmatrix}$, which correspond to the cases in which promoter $1$ transcribes repressor $2$ and in which promoter $2$ transcribes repressor $1$, respectively. The unstable equilibrium is $E_2=\begin{pmatrix} 0.8807 & 0.7808 \end{pmatrix}$. The separation between the DOAs corresponding to $E_1$ and $E_3$ is achieved by the stability boundary or separatrix. It is known that the unstable equilibrium $E_2$ belongs to the stability boundary \cite{Chiang88}. 

The problem of computing the DOAs of the two stable equilibria was previosuly addressed by means of PWA approximating dynamics in \cite{Yordanov2013} by computing reachable sets via a linear temporal logic formalism, however the method is very computationally expensive even for such a system with two states. 
By following the procedure in Section~\ref{sec:04.02}, the DOAs corresponding to $E_1$ and $E_3$ where computed and are shown in Figure~\ref{fig:05.03}.
Note that the computation procedure needs to be carried out for each system resulting by translating each nonzero equilibrium to the origin. For each system corresponding to $E_1$ and $E_3$, a quadratic FTLF candidate was considered, where the matrix $P$ is the solution of the classical Lyapunov inequality $A^\top P+PA<-I$, where $I$ denotes the identity matrix and
$A= \left[\frac{\partial f(x)}{\partial x}\right]_{x=E_i}$ for each $i=1,3$.
The blue level set defined by $C_1=0.07$ and $d=1.2$ corresponds to $E_1$ and the black level set defined by $C_3=0.8$ and $d=0.4$ corresponds to $E_3$. The trajectories starting from initial conditions close to the stability boundary will go to $E_2$ and from there via its unstable directions, they will converge  either to $E_1$ or $E_3$. As shown in Figure~\ref{fig:05.03}, the DOA estimates go very close to $E_2$ (red), thus close to the stability boundary. However, the computed sets seem conservative with respect to the directions of the vector fields for initial conditions far from the equilibria in the positive orthant. This is due to the fact that higher level sets of the corresponding LFs would intersect with the unstable equilibrium and violate the stability boundary. In fact, this is not an issue for the toggle switch as the real--life behavior is centered around the three equilibria and stability boundary.
\subsection{Nonpolynomial 3D system--the HPA axis}
\label{ex:05.04}
The following model has been proposed in \cite{Andersen2013} to illustrate the behavior of the  Hypothalamic-Pituitary-Adrenal (HPA) axis.
\begin{eqnarray}
\label{eq:05.04}
\nonumber
  \dot{x}_1  &=& \left(1+\xi\frac{x_3^\alpha}{1+x_3^\alpha}
  -\psi\frac{x_3^\gamma}{x_3^\gamma +\tilde{c}_3^\gamma}\right)
  -\tilde{w}_1x_1 \\
  \dot{x}_2  &=&\left(1-\rho\frac{x_3^\alpha}{1+x_3^\alpha}\right)x_1-\tilde{w}_2x_2\\
  \nonumber
  \dot{x}_3  &=& x_2-\tilde{w}_3x_3.
\end{eqnarray}
The HPA axis is a system which acts mainly at maintaining body homeostasis by regulating the level of cortisol. The three hormones involved in the HPA axis are the CRH ($x_1$), the ACTH ($x_2$) and the cortisol ($x_3$). For certain parameter values the system \eqref{eq:05.04} has a unique  stable equilibrium, which relates to cortisol level returning to normal after periods of mild stress in healthy individuals. When the parameters are perturbed, bifurcation can occur which leads to bistability. The stable states correspond to hypercortisolemic and hypocortisolemic  equilibria, respectively.
\begin{figure}[htp]
\centerline{\includegraphics[width=5.0cm,height=4.8cm]{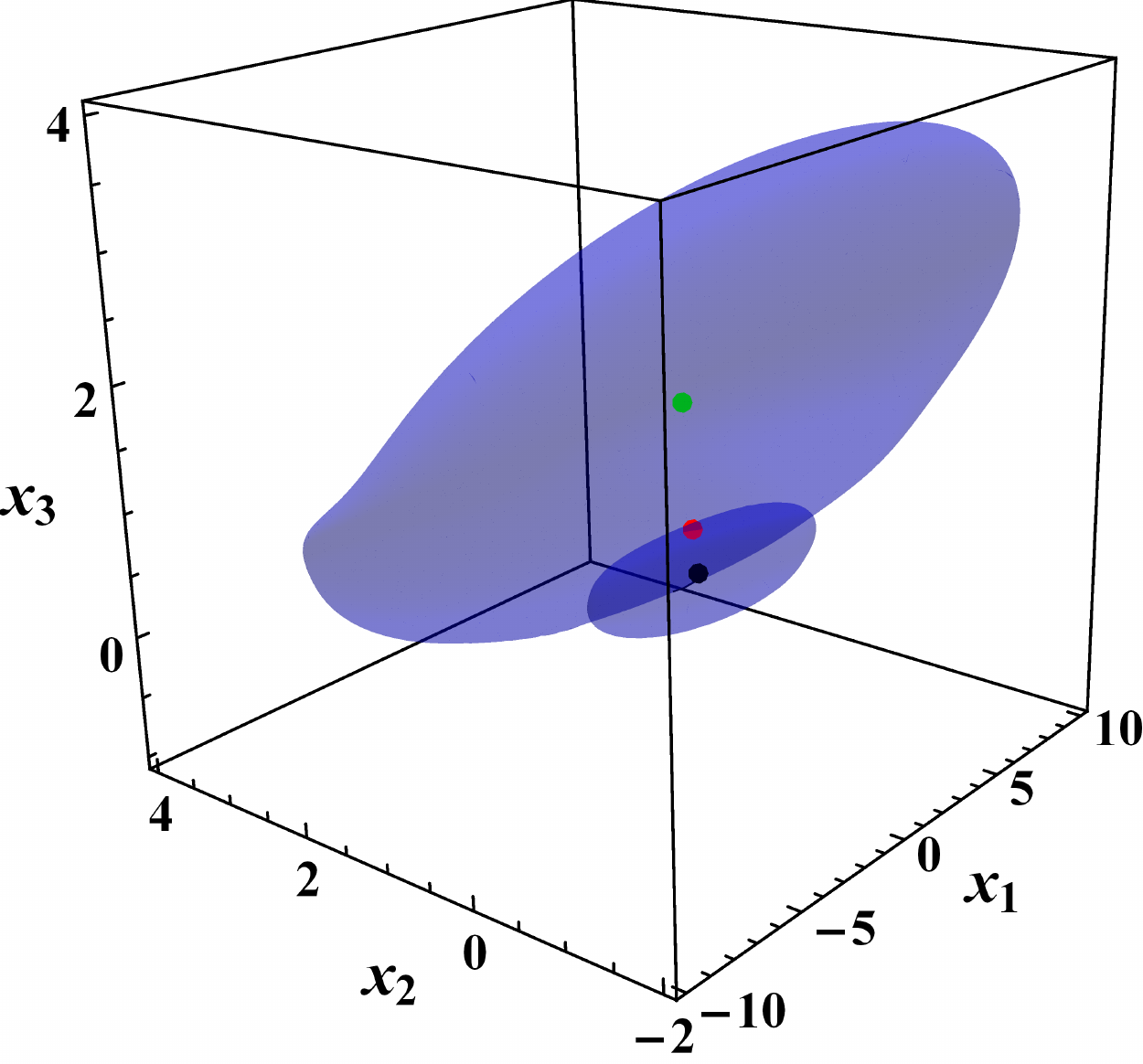}}
\caption{Level sets of the computed LFs corresponding to $E_1$--green and $E_3$--black for the HPA system \eqref{eq:05.04}. }
\label{fig:05.04}
\end{figure}
For the parameter values $\tilde{w}_1=4.79$, $\tilde{w}_2=0.964$, $\tilde{w}_3=0.251$, $\tilde{c}_3=0.464$, $\psi=1$, $\xi=1$, $\rho = 0.5$, $\gamma =\alpha=5$, the HPA system has three equilibria, $E_1 = \begin{pmatrix}0.1170 & 0.1199 & 0.4778\end{pmatrix}$,
$E_2 = \begin{pmatrix} 0.2224 & 0.2017 & 0.8039\end{pmatrix}$ and
$E_3 = \begin{pmatrix} 0.7833 & 0.4316 & 1.7196 \end{pmatrix}$, with $E_1$ and $E_3$ stable and $E_2$ unstable. By following the same steps as those indicated in the case of the toggle switch system, for the $E_1$ equilibrium the level set defined by $C_1=0.08$ and $d=0.4$ is plotted in Figure~\ref{fig:05.04} (the lower set) and for the $E_2$ equilibrium the level set defined by $C_2=1$ and $d=0.4$ is plotted in Figure~\ref{fig:05.04} (the upper set).

It is worth noting that for this system, the logarithmic norm defined by taking the $2$--norm, $\mu_2(A)$ is positive, with $A=\left[\frac{\partial f}{\partial x}\right]_{x=E_1}$. Thus, equality \eqref{eq:04.03} will not be satisfied. However, since the considered FTLF is $V(x)=x^\top P x$, then the weighted logarithmic norm satisfies $\mu_{2,P}(A)<0$.

\subsection{Nonpolynomial 3D system--the repressilator}
\label{ex:05.05}
Regulatory molecular networks, especially the oscillatory networks, have attracted a lot of interest from biologists and biophysiscists because they are found in many molecular pathways. 
Abnormalities of these processes lead
to various diseases, from sleep disorders to cancer. 
The naturally  occurring regulatory networks are very complex,
so their dynamics have been studied by highly
simplified models \cite{Buse2010}, \cite{Elowitz2000}. 
These models are particularly valuable because  they can provide an understanding of the important properties in the naturally occurring regulatory networks and, thus, support the engineering of artificial ones. Moreover, these rather simple, models can describe behaviors observed in experiments rather well \cite{Elowitz2000}.

An example of such a network is the repressilator.  Its genetic implementation uses
three proteins that cyclically repress the synthesis of one another. 
\begin{figure}[!h]
\centering
\subfigure[The level set $W(x)=0.32$--blue and $\dot{W}=0$--red.]{\includegraphics[width=0.42\columnwidth]{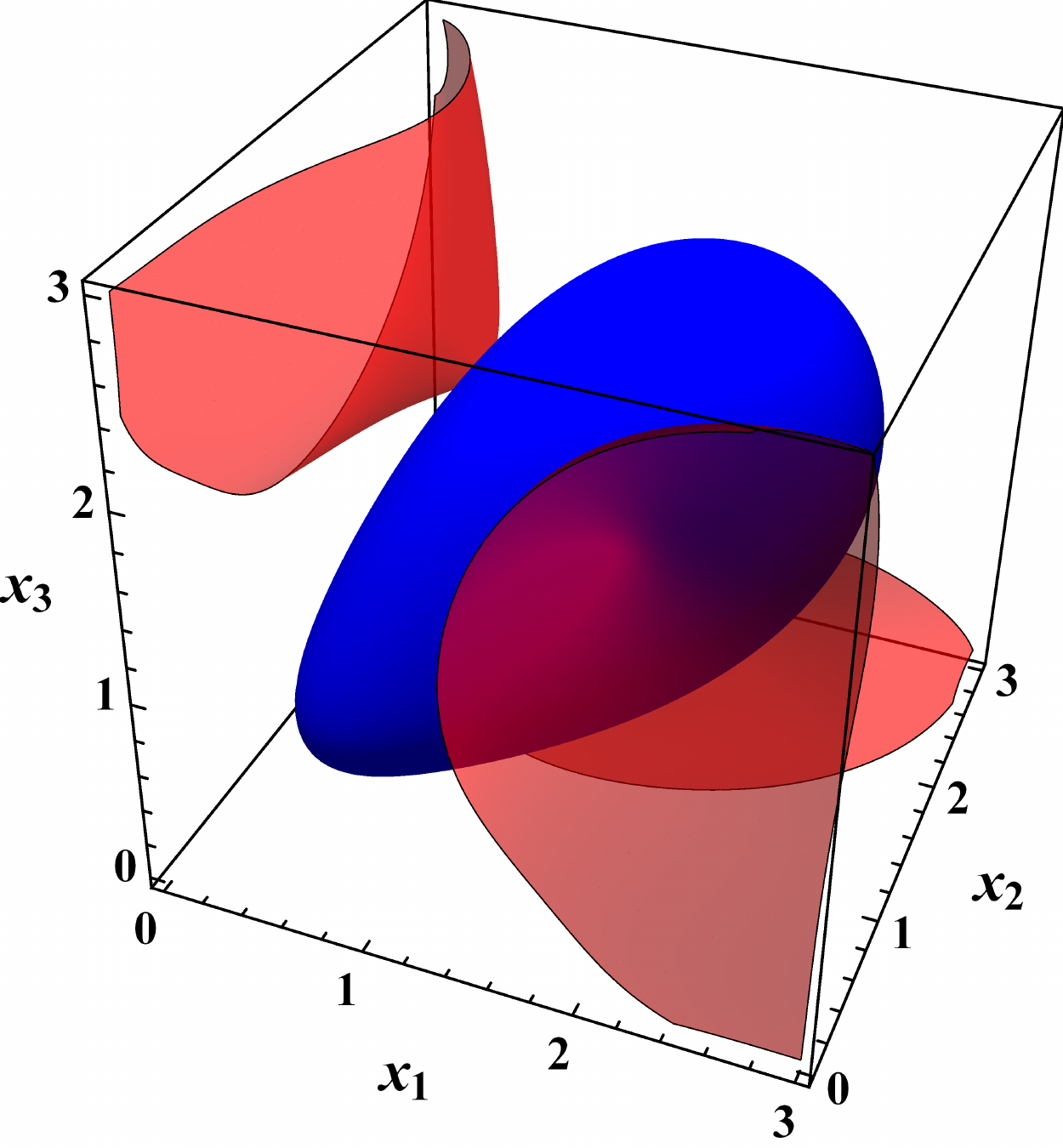}
\label{ch6:fig:03.01}}
\hspace{0.1cm}
\subfigure[The level set $W(x)=0.32$ and some trajectories of \eqref{ch6:eq:03.01}.]{\includegraphics[width=0.42\columnwidth]{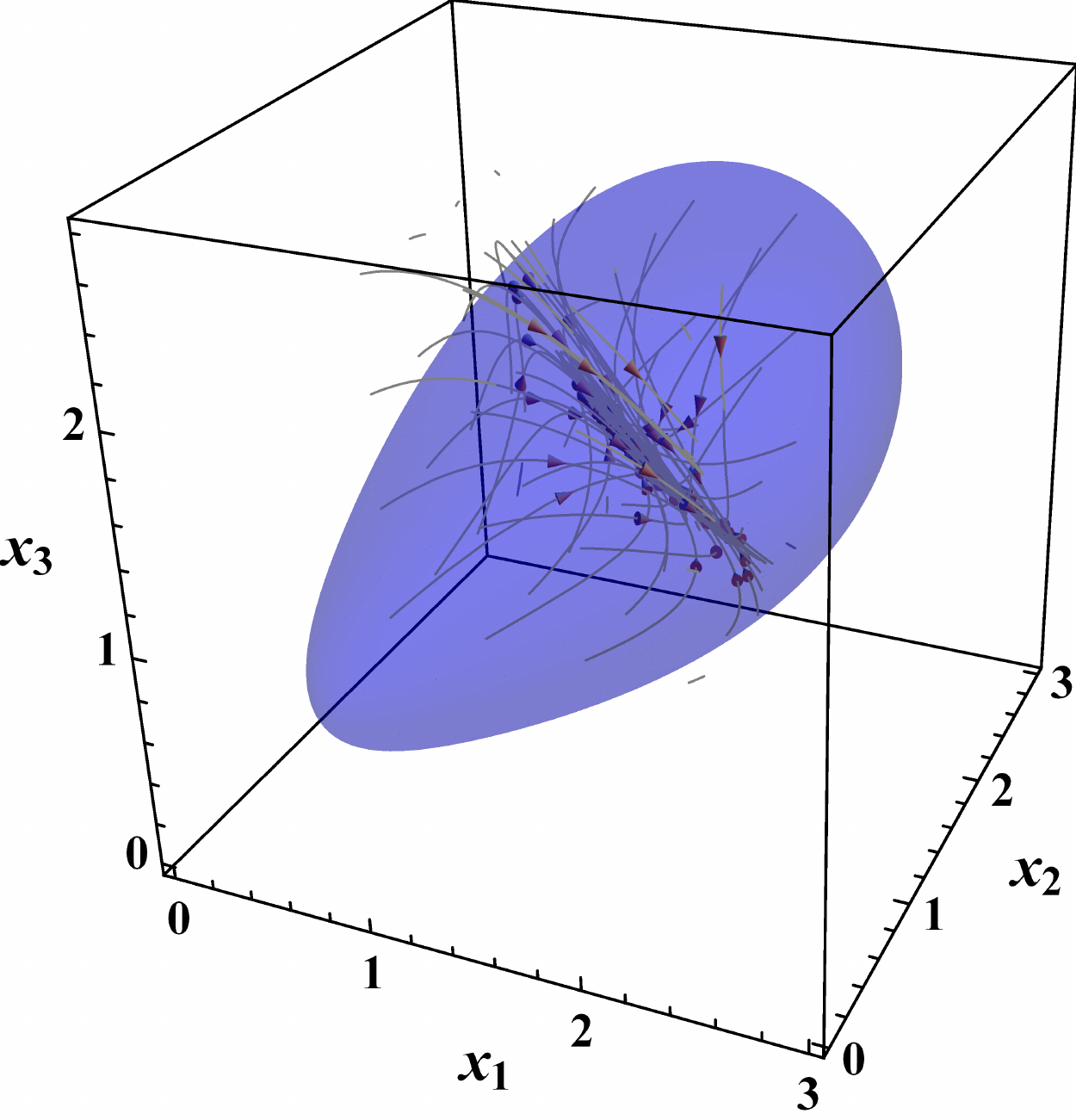}
\label{ch6:fig:03.02}}
\caption{Computation of the DOA of $E=(1.516, 1.516, 1.516)$ for \eqref{ch6:eq:03.01}.}
\end{figure} 
A model for the reprissilator was first proposed in \cite{Elowitz2000} and we consider here the simplified version from \cite{Buse2010}:
\begin{eqnarray}
\label{ch6:eq:03.01}
\nonumber
  \dot{x}_1  &=&\frac{\alpha}{1+x_2^\beta}-x_1  \\
  \dot{x}_2  &=&\frac{\alpha}{1+x_3^\beta}-x_2 \\
  \nonumber
  \dot{x}_3  &=& \frac{\alpha}{1+x_1^\beta}-x_3 .
\end{eqnarray}
The states $x_1$, $x_2$, and $x_3$ are proportional to protein concentrations. All negative terms in
the right-hand side represent degradation of the molecules.
The nonlinear function $h(x)=\frac{1}{1+x^\beta}$ reflects synthesis of the
mRNAs from the DNA controlled by regulatory elements
called promoters. $\beta$ is called cooperativity and reflects multimerization of the protein required to affect the promoter.
The three proteins are assumed to be identical, rendering the model symmetric and the order in choosing the states $x_1$, $x_2$ and $x_3$ does not influence the analysis outcome. 

In \cite{Buse2010} it was shown that for $\alpha>0$ and $\beta>1$, there is only one equilibrium point for the system \eqref{ch6:eq:03.01}, of the type $E=(r,r,r)$, where $r$ satisfies the equation $r^{\beta+1} +r-\alpha=0$. 

For the values $\alpha=5$ and $\beta =2$, $r=1.516$ and the eigenvalues of its corresponding linearized matrix are $\lambda=( -2.3936,
  -0.3032 + 1.2069i, -0.3032 - 1.2069i)$. By following the procedure in Section~\ref{sec:04.02} with the FTLF candidate $V(x)=x^\top P x$, with $P$ the identity matrix the value $d=0.4$ was obtained. For the obtained function $W=\int_0^d V(x+\tau f(x))\d\tau$ the level set given by $C=0.32$ was checked to be a subset of the true DOA of the system and it is plotted in Figure~\ref{ch6:fig:03.01} with blue. The zero level set of $\dot{W}=\nabla W^\top f(x)$, where $f(x)$ denotes the map describing \eqref{ch6:eq:03.01} is plotted in Figure~\ref{ch6:fig:03.01} with red. In Figure~\ref{ch6:fig:03.02} the level set $W(x)=0.32$ is shown together with some trajectories of \eqref{ch6:eq:03.01}. 
\subsection{Trigonometricc nonlinearity--the whirling pendulum}
\label{ex:05.06}
\begin{figure}[t!]
\centering
\includegraphics[width=0.6\columnwidth]{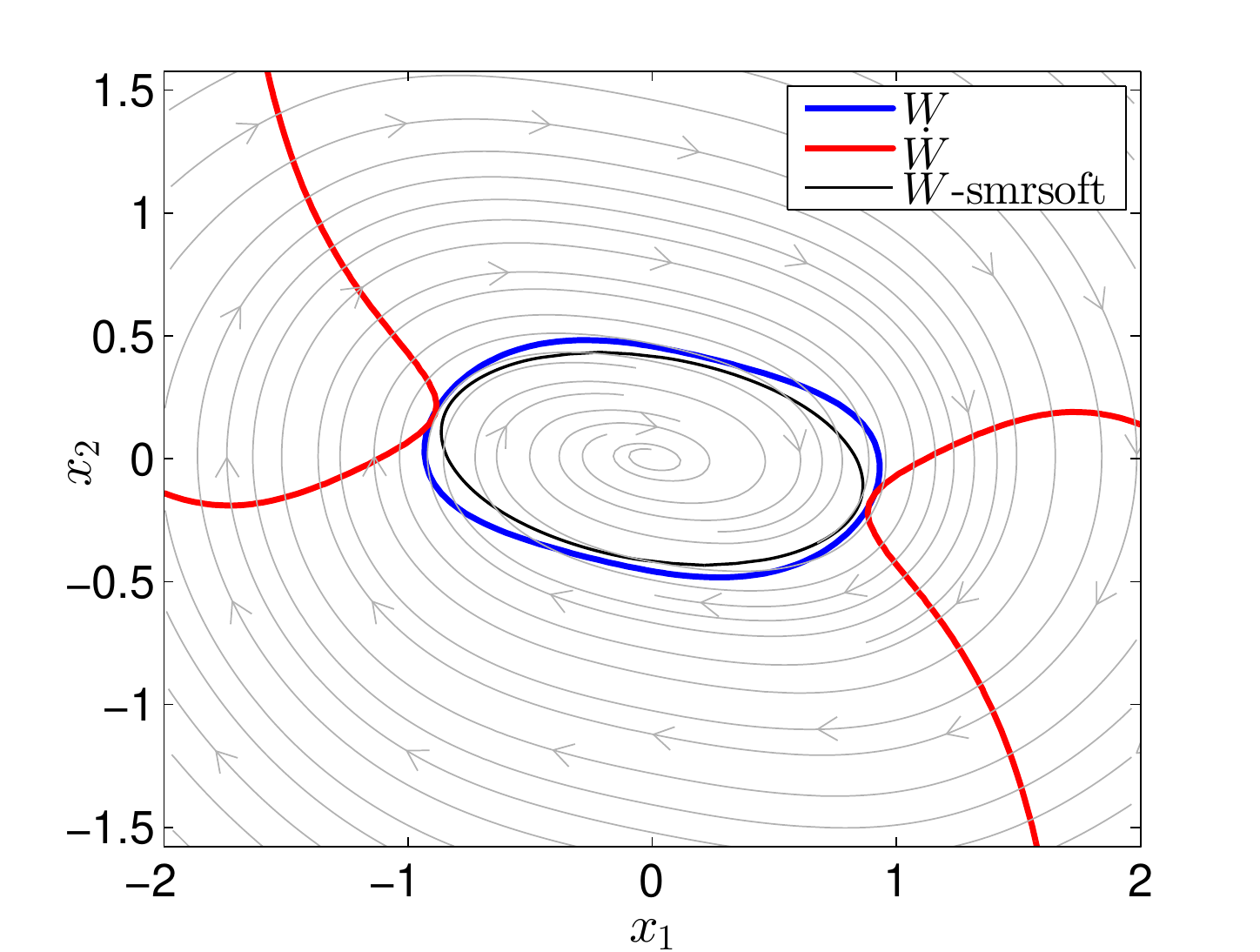}
\label{fig:05.06}
\caption{Level set of $W$ for $C=3.55$ computed for \eqref{eq:05.06} with  $d=1.1$,  plotted in blue, its corresponding derivative--red and the level set computed in \cite{Chesi2009} with the toolbox smrsoft. }
\end{figure}

We consider the system below, which was studied in \cite{Chesi2009} with the purpose to compute the DOA of the zero equilibrium.
\begin{eqnarray}
\label{eq:05.06}
\dot{x}_1 &=& x_2\\
\dot{x}_2 &=& \frac{-k_f}{mb} x_2 +\omega^2\sin(x_1)\cos(x_1) -\frac{g}{l_p}\sin(x_1),
\nonumber
\end{eqnarray}
  where $x_1$ is the angle with the vertical, $k_f = 0.2$ is the friction, $m_b = 1$ is the  mass of the rigid arm , $l_p=10$ is the length of the rigid arm, $\omega =0.9$ is the angular velocity and  $g=10$ is the gravity acceleration.
Therein a polynomial LF was computed, whose level set rendering a DOA estimate is shown in Figure~\ref{fig:05.06} with the black contour.
Following the same steps as in the previous examples, a Lyapunov function $W$ was computed on basis of the FTLF $V(x)=x^\top P x$, with
$P= \begin{pmatrix}
   3.6831  &  2.3169\\
    2.3169  & 14.7694
\end{pmatrix}$ and $d=1.1$.
The level set $C=3.55$ of $W(x)$ defines an estimate of the true DOA of the origin equilibrium of \eqref{eq:05.06} and it is shown with blue in Figure~\ref{fig:05.06} together with a vector field plot of the system.

\subsection{Trigonometric nonlinearity--multiple equilibria}
\label{ex:05.07}
Consider the system:
\begin{equation}
\label{eq:05.07}
\begin{split}
\dot{x_1}&=x_2\\
\nonumber
\dot{x_2}&=0.301 - \sin(x_1+0.4136) +\\
&\quad\,\,0.138\sin 2(x_1+0.4136)-0.279 x_2,
\end{split}
\end{equation}
\begin{figure}[htp]
\centering
\includegraphics[width=0.65\columnwidth]{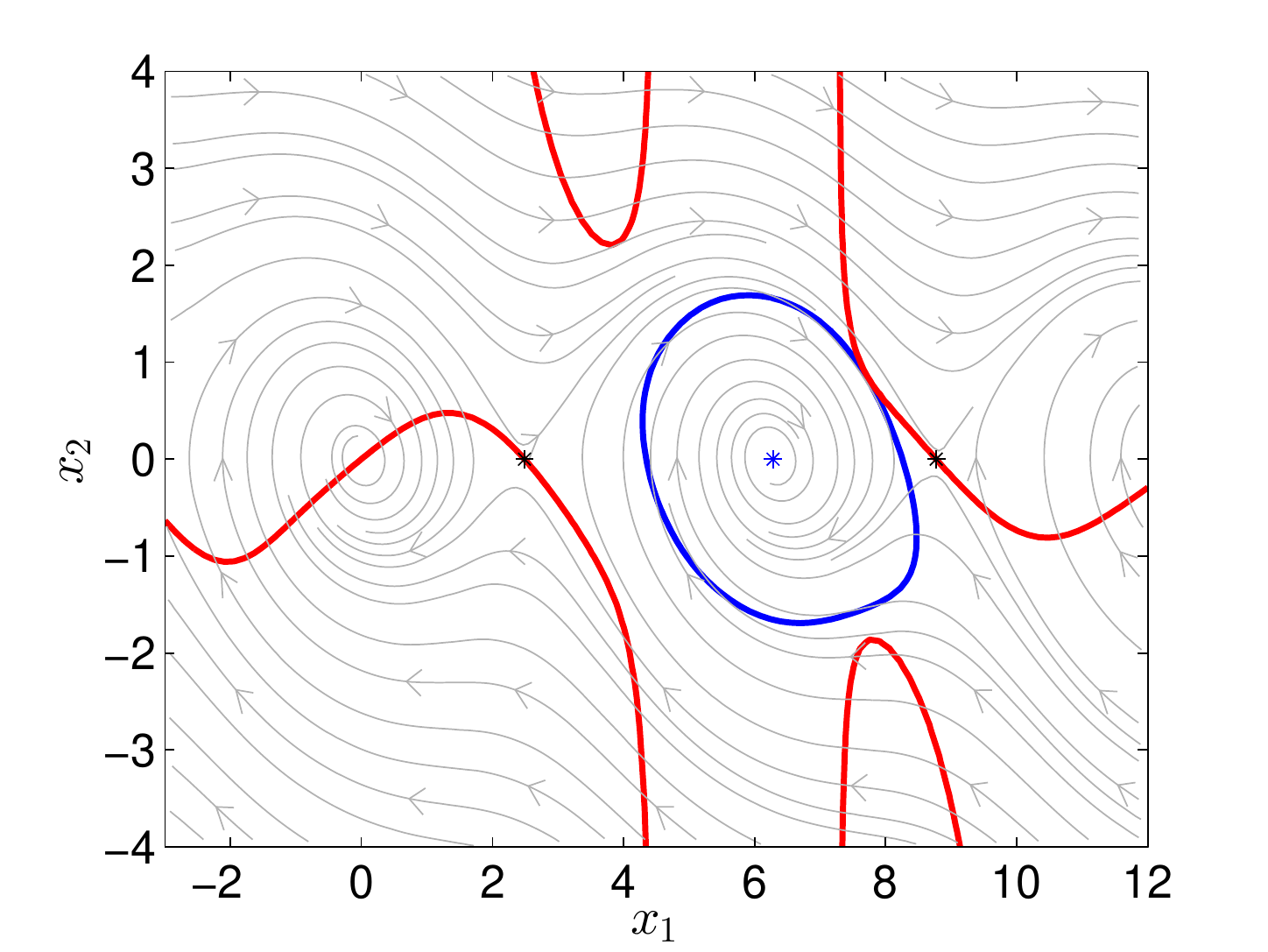}
\caption{The level set  $W(x)=5$--blue, its derivative $\dot{W}=0$--red, the stable equilibrium $E_1$--blue and the unstable ones $E_1, E_2$--black together with the vector field plot for \eqref{eq:05.07}.}
\label{fig:05.07}
\end{figure}
with the stable equilibrium $E_1=(6.284098\quad 0)^\top$, and the unstable ones $E_2=(2.488345 \quad 0)^\top$, $E_2=(8.772443 \quad 0)^\top$.  This system  was studied in \cite{Chiang88} with the purpose to compute the stability boundary of $E_1$.  By applying the steps described in Section~\ref{sec:04.02} for verification starting with the linearized dynamics,  for $V(x)=x^\top P x$ where $P=\begin{pmatrix}
    1.6448  &  0.3430\\
    0.3430  &  2.1255
\end{pmatrix}$, 
the condition \eqref{eq:norm} holds with $d=0.8$. The resulting LF $W=\int_{0}^d V(x+\tau f(x))\d\tau$ leads to the DOA estimate defined by $C= 5$ and plotted in Figure~\ref{fig:05.06} with blue contour. The zero level set of the corresponding derivative $\dot{W}=\nabla W(x)^\top f(x)=0$ is shown with red. 
\section{Conclusions}
\label{sec:06}
We have provided a new Massera type of LF which is enabled by imposing a finite--time condition on an arbitrary candidate function defined by a $\mc{K}_\infty$--function of the state norm. As the finite--time condition is verifiable numerically, this definition of the LF allows for an implementable algorithm towards obtaining nonconservative DOA estimates.
Compared to classical constructions which require for finite time interval integration exponential stability or they allow for $\mc{K}\mc{L}$--stability by imposing a specific function under the integral,
we provide an approach which provides two more degrees of freedom. On one hand, we allow for $\mc{K}\mc{L}$--stability (under Assumption~\ref{as:02.01}) and on the other hand the construction of the LF is based on any $\mc{K}_\infty$ function of the norm of the solution of the system. For future work we are considering building a Matlab toolbox for computing LFs and DOA estimations for nonlinear systems. 




\bibliographystyle{IEEEtran}        
\bibliography{AlinaRefs,AlinaRefsCh4}           

\begin{thebibliography}{10}
\providecommand{\url}[1]{#1}
\csname url@rmstyle\endcsname
\providecommand{\newblock}{\relax}
\providecommand{\bibinfo}[2]{#2}
\providecommand\BIBentrySTDinterwordspacing{\spaceskip=0pt\relax}
\providecommand\BIBentryALTinterwordstretchfactor{4}
\providecommand\BIBentryALTinterwordspacing{\spaceskip=\fontdimen2\font plus
\BIBentryALTinterwordstretchfactor\fontdimen3\font minus
  \fontdimen4\font\relax}
\providecommand\BIBforeignlanguage[2]{{%
\expandafter\ifx\csname l@#1\endcsname\relax
\typeout{** WARNING: IEEEtran.bst: No hyphenation pattern has been}%
\typeout{** loaded for the language `#1'. Using the pattern for}%
\typeout{** the default language instead.}%
\else
\language=\csname l@#1\endcsname
\fi
#2}}

\bibitem{Massera1949}
J.~L. Massera, ``On {L}iapounoff's conditions of stability,'' \emph{Annals of
  Mathematics}, vol.~50, no.~3, pp. 705--721, 1949.

\bibitem{Kurzweil}
J.~Kurzweil, ``On the inversion of {L}japunov's second theorem on stability of
  motion,'' \emph{Czechoslovak Mathematical Journal}, vol.~81, pp. 455 -- 484,
  1956, {E}nglish transaltion in Americal Mathematical Society Transaltions
  (2), v.24, pp. 19-77.

\bibitem{Zubov64}
V.~I. Zubov, \emph{Methods of A. M. {L}yapunov and their application}.\hskip
  1em plus 0.5em minus 0.4em\relax Noordhoff, Groningen, 1964.

\bibitem{Hahn67}
W.~Hahn, \emph{Stability of motion}.\hskip 1em plus 0.5em minus 0.4em\relax Die
  Grundlehren der mathematischen Wissenschaften, Band 138, Springer, Berlin,
  1967.

\bibitem{VanVid85}
A.~Vannelli and M.~Vidyasagar, ``Maximal {L}yapunov functions and domains of
  attraction for autonomous nonlinear systems,'' \emph{Automatica}, vol.~21,
  no.~1, pp. 69--80, 1985.

\bibitem{Yoshizawa1966}
T.~Yoshizawa, \emph{Stability Theory by {L}iapunov's Second Method}.\hskip 1em
  plus 0.5em minus 0.4em\relax Publications of the Mathematical Society of
  Japan, No. 9, 1966.

\bibitem{Krasovskii1963}
N.~N. Krasovskii, \emph{Stability of motion: {A}pplications of {L}yapunov's
  second method to differential systems and equations with delay}.\hskip 1em
  plus 0.5em minus 0.4em\relax Stanford University Press, 1963.

\bibitem{LaSalleLef61}
J.~LaSalle and S.~Lefschetz, \emph{Stability by {L}iapunov's Direct Method with
  Applications}.\hskip 1em plus 0.5em minus 0.4em\relax Academic Press, New
  York, London, 1961.

\bibitem{Kalman1959}
R.~Kalman and J.~Bertram, ``Control system analysis and design via the ``second
  method'' of {L}yapunov: Part {I}, continuous-time systems,''
  \emph{Transactions of the AMSE, Series D: Journal of Basic Engineering},
  vol.~82, no.~2, pp. 371--393, 1960.

\bibitem{LinSontagWang1996}
Y.~Lin, E.~D. Sontag, and Y.~Wang, ``A smooth converse {L}yapunov theorem for
  robust stability,'' \emph{SIAM Journal on Control and Optimization}, vol.~34,
  no.~1, pp. 124--160, 1996.

\bibitem{Rosier92}
L.~Rosier, ``Homogeneous {L}yapunov function for homogeneous continuous vector
  field,'' \emph{Syst. Control Lett.}, vol.~19, no.~6, Dec. 1992.

\bibitem{TeelPraly}
A.~R. Teel and L.~Praly, ``A smooth {L}yapunov function from a
  class-$\mc{K}\mc{L}$ estimate involving two positive semidefinite
  functions,'' \emph{ESAIM: Control, Optimisation and Calculus of Variations},
  vol.~5, pp. 313--367, 2000.

\bibitem{KelletTeel2004}
C.~M. Kellett and A.~R. Teel, ``Weak converse {L}yapunov theorems and
  {C}ontrol-{L}yapunov functions,'' \emph{SIAM Journal on Control and
  Optimization}, vol.~42, no.~6, pp. 1934--1959, 2004.

\bibitem{Sontag1983}
E.~D. Sontag, ``A {L}yapunov-like characterization of asymptotic
  controllability,'' \emph{SIAM Journal on Control and Optimization}, vol.~21,
  no.~3, pp. 462--471, 1983.

\bibitem{Hafstein2007}
S.~F. Hafstein, ``An algorithm for constructing {L}yapunov functions,''
  \emph{Electronic Journal of Differential Equations, Monograph 08}, 2007.

\bibitem{Bjornsson2014}
J.~Bjornsson, P.~Giesl, S.~F. Hafstein, C.~M. Kellett, and H.~Li, ``Computation
  of continuous and piecewise affine {L}yapunov functions by numerical
  approximations of the {M}assera construction,'' in \emph{53rd IEEE Annual
  Conference on Decision and Control (CDC)}, Dec 2014, pp. 5506--5511.

\bibitem{Bjornsson2015}
J.~Bjornsson, P.~Giesl, S.~F. Hafstein, and C.~K.~M. Kellett, ``Computation of
  {L}yapunov functions for systems with multiple local attractors,''
  \emph{Discrete and Continuous Dynamical Systems}, vol.~35, no.~9, pp.
  4019--4039, 2015.

\bibitem{SiggiKelletLi2014}
S.~F. Hafstein, C.~M. Kellett, and H.~Li, ``Continuous and piecewise affine
  {L}yapunov functions using the {Y}oshizawa construction,'' in \emph{American
  Control Conference (ACC), 2014}, June 2014, pp. 548--553.

\bibitem{Hafstein2016}
------, ``Computing continuous and piecewise affine {L}yapunov functions for
  nonlinear systems,'' \url{https://epub.uni-bayreuth.de/1962/}, Bayreuth,
  2015.

\bibitem{Kellet2015survey}
C.~M. Kellett, ``Classical converse theorems in {L}yapunov's second method,''
  \emph{Discrete and Continuous Dynamical Systems - Series B}, vol.~20, no.~8,
  pp. 2333--2360, 2015.

\bibitem{Giesl2015}
P.~Giesl and S.~F. Hafstein, ``Review on computational methods for {L}yapunov
  functions,'' \emph{Discrete and Continuous Dynamical Systems - Series B},
  vol.~20, no.~8, pp. 2291--2331, 2015.

\bibitem{Khalil2002}
H.~K. Khalil, \emph{Nonlinear Systems}.\hskip 1em plus 0.5em minus 0.4em\relax
  Prentice Hall, 2002.

\bibitem{Aeyels1998}
D.~Aeyels and J.~Peuteman, ``A new asymptotic stability criterion for nonlinear
  time-variant differential equations,'' \emph{IEEE Transactions on Automatic
  Control}, vol.~43, no.~7, pp. 968--971, Jul 1998.

\bibitem{Geiselhart2014}
R.~Geiselhart, R.~H. Gielen, M.~Lazar, and F.~R. Wirth, ``An alternative
  converse {L}yapunov theorem for discrete-time systems,'' \emph{Systems \&
  Control Letters}, vol.~70, pp. 49 -- 59, 2014.

\bibitem{Soderlind}
G.~S{\"o}derlind, ``The logarithmic norm. {H}istory and modern theory,''
  \emph{BIT Numerical Mathematics}, vol.~46, no.~3, pp. 631--652, 2006.

\bibitem{Hu2004}
G.-D. Hu and M.~Liu, ``The weighted logarithmic matrix norm and bounds of the
  matrix exponential,'' \emph{Linear Algebra and its Applications}, vol. 390,
  pp. 145 -- 154, 2004.

\bibitem{Sontag98}
E.~D. Sontag, ``Comments on integral variants of {ISS},'' \emph{Systems \&
  Control Letters}, vol.~34, no. 1-2, pp. 93 -- 100, 1998.

\bibitem{Sontag1998}
------, \emph{Mathematical Control Theory: Deterministic Finite Dimensional
  Systems. Second Edition}, ser. Textbooks in Applied Mathematics.\hskip 1em
  plus 0.5em minus 0.4em\relax Springer, New York, 1998.

\bibitem{Geiselhart2015}
R.~Geiselhart, ``Advances in the stability analysis of large-scale
  discrete-time systems,'' Ph.D. dissertation, The Julius Maximilians
  University of Wurzburg, 2015.

\bibitem{Browder1996}
A.~Browder, \emph{Mathematical analysis. An introduction}.\hskip 1em plus 0.5em
  minus 0.4em\relax Springer, 1996.

\bibitem{Karafyllis2012}
I.~Karafyllis, ``Can we prove stability by using a positive definite function
  with non sign-definite derivative?'' \emph{IMA Journal of Mathematical
  Control and Information}, vol.~29, no.~2, pp. 147--170, 2012.

\bibitem{Chiang89}
H.~D. Chiang and J.~S. Thorp, ``Stability regions of nonlinear autonomous
  dynamical systems: a constructive methodology,'' \emph{IEEE Transactions on
  Automatic Control}, vol.~34, pp. 1229--1241, 1989.

\bibitem{Dieudonne69}
J.~Dieudonn\'{e}, \emph{Treatise on Analysis, Vol. {I}, Foundations of Modern
  Analysis}.\hskip 1em plus 0.5em minus 0.4em\relax Academic Press, New York
  and London, 1969.

\bibitem{Sontag2014}
Z.~Aminzare and E.~D. Sontag, ``Contraction methods for nonlinear systems: A
  brief introduction and some open problems,'' in \emph{53rd IEEE Conference on
  Decision and Control (CDC)}, 2014, pp. 3835--3847.

\bibitem{Chesi2011}
G.~Chesi, \emph{Domain of Attraction: Analysis and Control via SOS
  Programming}, ser. Lecture Notes in Control and Information Sciences.\hskip
  1em plus 0.5em minus 0.4em\relax Springer, 2011.

\bibitem{Hachicho2007}
O.~Hachicho, ``A novel {LMI} -based optimization algorithm for the guaranteed
  estimation of the domain of attraction using rational {L}yapunov functions,''
  \emph{Journal of the Franklin Institute}, vol. 344, no.~5, pp. 535 -- 552,
  2007.

\bibitem{Parillo2011}
A.~A. Ahmadi, M.~Krstic, and P.~A. Parrilo, ``A globally asymptotically stable
  polynomial vector field with no polynomial {L}yapunov function,'' in
  \emph{50th IEEE Conference on Decision and Control and European Control
  Conference (CDC-ECC)}, Dec 2011, pp. 7579--7580.

\bibitem{SiggiMicnon2015}
J.~Bjornsson, S.~Gudmundsson, and S.~F. Hafstein, ``Class library in {C}++ to
  compute {L}yapunov functions for nonlinear systems,''
  \emph{IFAC-PapersOnLine}, vol.~48, no.~11, pp. 778 -- 783, 2015, 1st {IFAC}
  Conference on Modelling, Identification and Control of Nonlinear Systems,
  2015.

\bibitem{Gardner2000}
T.~S. Gardner, C.~R. Cantor, and J.~J. Collins, ``Construction of a genetic
  toggle switch in {E}scherichia coli,'' \emph{Nature}, vol. 403, no. 6767, pp.
  339--342, 2000.

\bibitem{Chiang88}
H.~D. Chiang, M.~W. Hirsch, and F.~W. Wu, ``Stability regions of nonlinear
  autonomous dynamical systems,'' \emph{IEEE Transactions on Automatic
  Control}, vol.~33, pp. 16--27, 1988.

\bibitem{Yordanov2013}
N.~Yordanov, J.~Tumova, I.~Cerna, J.~Barnat, and C.~Belta, ``Formal analysis of
  piecewise affine systems through formula-guided refinement,''
  \emph{Automatica}, vol.~49, no.~1, pp. 261 -- 266, 2013.

\bibitem{Andersen2013}
M.~Andersen, F.~Vinther, and J.~T. Ottesen, ``Mathematical modeling of the
  hypothalamic-pituitary-adrenal gland (hpa) axis, including hippocampal
  mechanisms,'' \emph{Mathematical Biosciences}, vol. 246, no.~1, pp. 122 --
  138, 2013.

\bibitem{Buse2010}
O.~Buse, R.~P\'erez, and A.~Kuznetsov, ``Dynamical properties of the
  repressilator model,'' \emph{Phys. Rev. E}, vol.~81, pp.
  066\,206--1--066\,206--7, June 2010.

\bibitem{Elowitz2000}
M.~B. Elowitz and S.~Leibler, ``A synthetic oscillatory network of
  transcriptional regulators,'' \emph{Nature}, vol. 403, no. 6767, pp.
  335--338, 2000.

\bibitem{Chesi2009}
G.~Chesi, ``Estimating the domain of attraction for non-polynomial systems via
  {LMI} optimizations,'' \emph{Automatica}, vol.~45, no.~6, pp. 1536 -- 1541,
  2009.

\end{thebibliography}
\end{document}